%% file: ARXIV_EFR_OPINF.tex
\documentclass[a4paper,11pt]{article}
\usepackage[numbers]{natbib}
 
\usepackage{amsfonts,amsmath,amssymb,amsthm}
\usepackage{pifont}
\usepackage{graphicx}
\usepackage[table]{xcolor}

\usepackage{bm}
\usepackage{url}
\usepackage{subcaption}
\usepackage{comment}
\usepackage[colorlinks=true,urlcolor=blue,linkcolor=blue,citecolor=blue]{hyperref}
\usepackage{booktabs}
\usepackage{mathtools}
\usepackage{algorithm}
\usepackage{algpseudocode}
\usepackage{combelow}
\usepackage[normalem]{ulem}

\newcommand{\xmark}{\ding{55}}%
 
\usepackage[capitalise,noabbrev]{cleveref}[]
 
\begin{document}

\title{Learning Long-Term Stable Operator Inference Reduced-Order Models of Fluid Flows through Online Spatial Filtering}

\input{Notation.tex}

\author{
    Ian Moore\thanks{Department of Computational Applied Mathematics and Operations Research, Rice University, Houston, TX 77005, USA (\texttt{ian.moore@rice.edu}).} \and
    Ping-Hsuan Tsai\thanks{Department of Applied Mathematics, National Yang Ming Chiao Tung University, Hsinchu 30010, Taiwan.} \and
    Anthony Gruber\thanks{Sandia National Laboratories, Albuquerque, NM 87123, USA.} \and
    Ionu\c{t} Farca\c{s}\thanks{Department of Mathematics and Division of Computational Modeling and Data Analytics, Virginia Tech, Blacksburg, VA 24061, USA.} \and
    Christopher Wentland\footnotemark[3] \and
    Irina Tezaur\footnotemark[3] \and
    Traian Iliescu\footnotemark[4]
}
 
\maketitle
 
\begin{abstract}
This paper introduces an online evolve--filter--relax (EFR) strategy for long-term stability of Operator Inference (OpInf) reduced-order models (ROMs) of complex fluid flow simulations.
The main novelty of the new EFR-OpInf strategy is the use of online (i.e., at the learned ROM online evaluation level) spatial filtering inspired from large eddy simulation to significantly improve long-term stability and predictive performance of standard OpInf.
Furthermore, the EFR-OpInf strategy reduces, and in some cases even eliminates, the need for standard $L^2$ regularization, while providing physical interpretability for the OpInf hyperparameters.
The EFR-OpInf framework is modular, readily integrated into existing OpInf workflows, and accommodates a user-selected ROM filtering strategy.
We demonstrate the EFR-OpInf's effectiveness using a fully non-intrusive projection-based ROM filter, a ROM differential filter, and a hybrid projection-differential ROM filter.
The new EFR-OpInf models are evaluated on a high-P\'eclet-number convection--diffusion--reaction problem that embeds the variation in one parameter and two unsteady Navier--Stokes problems that focus on predictions beyond a training horizon:
a transitional two-dimensional flow past a cylinder and a three-dimensional turbulent minimal channel flow.
Across these three scenarios, EFR-OpInf can reduce prediction errors by up to an order of magnitude relative to standard OpInf.
Moreover, EFR-OpInf remains stable over long prediction horizons in cases where standard OpInf diverges.
Depending on the ROM filter used, EFR-OpInf's computational cost is comparable to that of standard OpInf.
\end{abstract}

\section{Introduction}
High-fidelity numerical simulations of fluid flows often require the solution of large-scale nonlinear systems over many timesteps, making them computationally expensive even on modern computing architectures. 
This cost becomes particularly restrictive in the many-query settings common to uncertainty quantification, optimization, parameter studies, inverse problems, and real-time prediction and control, where the full-order model (FOM) must be evaluated many times or under strict computational constraints. 
Reduced-order models (ROMs) address this challenge by replacing the high-dimensional FOM with a low-dimensional approximation that retains dominant flow features while substantially reducing the computational cost. 

Among the most established approaches, standard intrusive Galerkin ROMs (G-ROMs)~\cite{Benner2015ROM,Benner2021ROM,hesthaven2015certified,quarteroni2015reduced} construct such low-dimensional approximations in two stages:
an offline stage that projects high-dimensional operators onto a low-dimensional space, and an online stage that evolves the resulting low-dimensional model in time.  
Due to its use of reduced governing equations, this can yield a substantial speed-up over the FOM numerical solution with a minimal loss of accuracy. 
While G-ROMs have been successfully used in the numerical simulation of laminar fluid flows~\cite{holmes2012turbulence,noack2011reduced,rozza2022advanced}, accurately representing these complex dynamics 
generally requires a prohibitively large number of ROM basis functions. 
Retaining computational efficiency in the ROM therefore requires the use of a relatively low-dimensional approximation; however, this 
does not adequately capture the effects of the discarded modes on the resolved dynamics. 
This loss of information can lead to a non-physical accumulation of energy and, consequently, spurious numerical oscillations in the reduced model.
Thus, for realistic, convection-dominated flows, G-ROMs are generally used in conjunction with closure models~\cite{ahmed2021closures} and/or stabilizations~\cite{BALAJEWICZ2016224,Parish2025residual} to account for the effect of the discarded ROM basis functions.

One of the main challenges in applying G-ROMs to realistic engineering, scientific, and biomedical flow problems is their intrusive nature. 
Intrusive 
models require access to the underlying FOM operators, which is challenging or even impossible in legacy, highly specialized, large-scale, or proprietary simulation codes. 
One of the most promising strategies to address this shortcoming of G-ROMs is {\it Operator Inference (OpInf)}~\cite{PEHERSTORFER2016_Opinf}, a scientific machine learning technique 
that bridges the gap between black-box machine learning and G-ROMs. 
Specifically, OpInf mimics the two-stage structure of G-ROMs but changes the offline stage from intrusive to non-intrusive: instead of using the FOM operators to assemble the ROM operators, OpInf learns reduced operators 
{through a non-intrusive least-squares minimization problem. 
The OpInf framework has enabled the use of projection-like ROMs for problems that were prohibitively challenging for intrusive methods, such as multiphysics combustion processes~\cite{McQuarrie2021Regularized, QFW21} and rotating detonation rocket engine simulations with more than $4 \times 10^{7}$ spatial degrees of freedom (DoFs)~\cite{farcas2023parametric,FARCAS2025109619}, where access to the FOM codes can be difficult or even impossible.

Despite the increased utility that comes with learning operators, there is a cost.  For complex applications, standard OpInf workflows require careful regularization of the operator-learning regression problem~\cite{farcas2023parametric, GKIMISKIS2025localizedopinf,McQuarrie2021Regularized, QFW21,tezaur2026hybridcouplingoperatorinference}, typically through tuning one or more Tikhonov hyperparameters.
The performance of the resulting OpInf model depends critically on these hyperparameters.  For example,  
the numerical investigation in Section~\ref{sec:numerical_results} shows that, for short training horizons, choosing hyperparameters that ensure the long-term stability of the resulting OpInf model
can require excessively large amounts of regularization, 
yielding inaccurate predictions. %
On the other hand, choosing hyperparameters that only minimize training error can yield OpInf models that overfit, i.e., ROMs that accurately reproduce the training data but become inaccurate and unstable during long-term prediction.
In addition to their critical effect on OpInf's stability and accuracy, these regularization parameters can span several orders of magnitude \cite{FARCAS2025109619,GKIMISKIS2025localizedopinf,tezaur2026hybridcouplingoperatorinference}, making OpInf hyperparameter selection challenging and computationally demanding (see, however, \cite{FARCAS2025109619} for an efficient, parallel implementation of the OpInf workflow including the hyperparameter selection).
Furthermore, these hyperparameters do not have a direct physical interpretation with respect to the ensuing ROM solution.

To address these challenges, we propose incorporating spatial filtering online into the OpInf solution to suppress spurious oscillations and promote stable predictions, potentially over long time horizons. 
Spatial filtering is a large eddy simulation- (LES-) inspired stabilization approach that is well established at the FOM level~\cite{berselli2006LES,Layton2012AD,Sagaut2006LES} and has also proved effective in intrusive G-ROMs~\cite{Moore2025_ADLROM,Wells2017EF,xie2018-NA_Leray}. 
In convection-dominated flow problems, for example, ROM spatial filtering has been used to mitigate the excessive energy accumulation caused by spatial under-resolution~\cite{Wells2017EF} and to suppress numerical oscillations that can lead to model blow-up during time integration~\cite{SANFILIPPO2023_ADL}. 
These methods provide a physically motivated alternative to purely algebraic stabilization techniques.  
The objective of this paper is to demonstrate that spatial filtering applied during the online stage can similarly enhance the long-term stability and predictive performance of OpInf ROMs of fluid flows.

To accomplish this, we introduce a new online filtering strategy: the {\it evolve-filter-relax OpInf (EFR-OpInf)} ROM. 
EFR is a three-step process that has successfully been applied at both the FOM~\cite{Bertagna2016_EFR,neda2012_EFR_analysis, strazzullo2021_EFRConsistency} and ROM levels~\cite{ GIRFOGLIO2023FILTER, IVAGNES2026_EFR,Wells2017EF}.   
EFR is a simple and modular strategy that takes an intermediate solution 
produced by a time-stepping scheme, filters this solution using a spatial filter of the user's choice, and then relaxes the filtered solution onto the tentative solution to provide the model update. 
This flexible strategy allows stabilization with physical meaning, filters with known parameter ranges, and numerical analysis that supports initial parameter choices. 
As shown in the numerical investigation in Section~\ref{sec:numerical_results}, the proposed EFR-OpInf model improves the performance of the standard OpInf in cases of limited training data availability. %
Indeed, for three challenging scenarios across two- and three-dimensional domains, EFR-OpInf: 
\begin{enumerate}
    \item[(i)]  %
        remains stable over long prediction horizons in cases for which standard OpInf diverges;    
    \item[(ii)] can reduce the standard OpInf prediction errors by up to an order of magnitude;
    \item[(iii)] can reduce, and in some cases even eliminate, the need for $L^2$ regularization in the OpInf process, %
    while introducing stabilization %
    filtering parameters with a clear physical interpretation; %
    \item[(iv)] %
    introduces only a modest computational overhead compared with standard OpInf; and  %
    \item[(v)] is a modular strategy that can be easily integrated into existing OpInf workflows. %
\end{enumerate}

The main novelty of our approach is the development of an online, ROM-level spatial filtering strategy for OpInf models. 
To %
our knowledge, this is the first use of online spatial filtering at the ROM level in the OpInf setting. 
The only previous work combining spatial filtering and OpInf was done in \cite{farcas2023filtering}, in which Gaussian filtering was applied to the FOM training data as an offline preprocessing step before constructing the OpInf ROM.
The proposed EFR-OpInf strategy is fundamentally different from 
this approach: 
in EFR-OpInf, spatial filtering is applied online to the evolving reduced-order solution rather than offline 
on the FOM data used for training. 
This distinction has two important practical consequences. 
First, the ROM-level filtering employed by EFR-OpInf scales with the reduced dimension, typically of order $\mathcal{O}(10)$, whereas FOM-level filtering used in \cite{farcas2023filtering} scales with the very large dimension of the full-order state.
Second, the online formulation allows filtering to be activated, deactivated, or adjusted during time integration, providing additional flexibility for balancing the accuracy and stability of OpInf predictions.

The remainder of this paper is organized as follows.
Section~\ref{sec:background} introduces notation and provides background on full-order modeling (Section~\ref{subsec:FOM}), OpInf with standard $L^2$ regularization (Section~\ref{subsec:OpInf_summary}), as well as the three filtering algorithms that we will use (Section \ref{subsec:filtering_methods}).
Section~\ref{sec:EFR_OpInf} presents the main methodological contribution of this work: the EFR-OpInf algorithm. 
Section~\ref{sec:numerical_results} then evaluates its performance against standard OpInf across three scenarios. 
We first consider a two-dimensional (2D) convection-diffusion-reaction (CDR) problem, 
constructing a parametric OpInf model 
used to predict a solution with a sharp boundary layer at P\'eclet number 6,000. 
We then investigate EFR-OpInf for the flow past a cylinder (FPC) arising from the 2D unsteady Navier--Stokes equations (NSE), demonstrating that the proposed model can produce stable future-state predictions from transitional flow data. 
We further assess the framework on a more challenging NSE scenario: the three-dimensional (3D) minimal channel flow  (MCF) in a turbulent regime. 
Overall, these results show that EFR stabilization can improve the long-time stability of OpInf models across increasingly complex flow settings.
Section~\ref{sec:conclusions} concludes the paper.

\section{Background} \label{sec:background}

This section introduces notation and provides background on full-order fluid flow modeling (Section~\ref{subsec:FOM}), OpInf with standard $L^2$ regularization (Section~\ref{subsec:OpInf_summary}), and the three filtering algorithms considered in the present work (Section~\ref{subsec:filtering_methods}).

\subsection{High-fidelity Fluid Flow Modeling} 
    \label{subsec:FOM}

Since the filtering strategies proposed in this work are motivated by LES for convection-dominated flows~\cite{berselli2006LES,Sagaut2006LES}, in the numerical investigation in Section~\ref{sec:numerical_results} we consider two convection-dominated problems: a high-P\'eclet number convection-diffusion-reaction equation, which is a simplified computational setting, and the convection-dominated NSE.
Since in the next sections the proposed filtering strategies will be developed for the NSE, in this section we briefly outline these equations. 
The unsteady incompressible NSE are written as follows:
\begin{equation}
\begin{split}
\frac{\partial \bu}{\partial t} - \mathrm{Re}^{-1} \Delta \bu  + (\bu \cdot \nabla) \bu + \nabla p &= {\boldsymbol{f}}, \quad\text{ in } \Omega \times (0, t_f] \\
\nabla \cdot \bu &= 0, \quad\text{ in } \Omega \times (0, t_f].
\label{eq:strong_NSE}
\end{split}
\end{equation}
The NSE~\eqref{eq:strong_NSE} are posed on a spatial domain $\Omega \in \mathbb{R}^d, d = 2 \text{ or } d = 3$ with velocity field $\bu: \Omega \times [0, t_f] \to \mathbb{R}^d$, where $t_f$ denotes the final time, %
pressure field $p:  \Omega \times [0, t_f] \to \mathbb{R}$, and forcing $\boldsymbol{f}$. 
We consider the continuous pressure and velocity spaces, 
$
    Q = L_0^2(\Omega) := \left\{ q \in L^2(\Omega) , \int_\Omega q = 0 \right\}, 
    \ 
    \mathbf{X} = H_0^1(\Omega;\ \mathbb{R}^d) := \left\{ \bu \in H^1(\Omega;\ \mathbb{R}^d), \bu|_{\partial\Omega}  = \mathbf{0} \right\}
$, 
suitable for %
the NSE~\eqref{eq:strong_NSE} equipped with homogeneous Dirichlet boundary conditions.  
The case of inhomogeneous Dirichlet conditions can be handled with minor adjustment to these definitions. 
For the %
2D flow past a cylinder test problem in Section~\ref{sec:NSE}, the %
FOM data are generated by solving the NSE %
with %
conforming finite element (FE) spaces, %
here consisting of
linear-quadratic Taylor--Hood element pairs. 
The spatial and temporal discretization parameters are provided in Section~\ref{sec:NSE}.
For the %
3D minimal channel flow in Section~\ref{sec:NSE_mfu}, the FOM data are generated using the spectral element method (SEM) to solve the NSE. 
The velocity and pressure are approximated using the $\mathbb{P}_N$--$\mathbb{P}_{N-2}$ formulation, and time integration is performed using the third-order semi-implicit BDF3/EXT3 scheme 
that treats the viscous term implicitly and the nonlinear term explicitly~\cite{fischer2017recent}. 
The spatial and temporal discretization parameters are provided in Section~\ref{sec:NSE_mfu}.
In what follows, $(\cdot,\cdot)$ denotes the $L^2$ inner product, while all other norms are specified explicitly when introduced.

\subsection{Operator Inference}
    \label{subsec:OpInf_summary}

OpInf is a scientific machine learning framework for constructing 
projection-like ROMs from data, typically used when the governing equations exhibit polynomial structure~\cite{Kramer2024_Opinf,PEHERSTORFER2016_Opinf}. 
For systems involving more general nonlinearities, lifting transformations can be used to expose polynomial structure in the transformed state variables, potentially through the introduction of auxiliary variables~\cite{QIAN20201LiftandLearn}. 
In the context of fluid flows governed by the NSE, OpInf seeks to learn 
quadratic ROMs from data generated by the semi-discrete governing equations~\eqref{eq:strong_NSE} (discretized via, e.g., FEM or SEM), thereby retaining the same polynomial degree as intrusive G-ROMs.

\subsubsection{Reduced basis and model form}

Let $[t_i,t_f]$ denote the time interval of interest, where $t_i > 0$ and $t_f$ denote the initial and final times, respectively. 
Moreover, let $N\in\mathbb{N}$ denote the dimension of the spatial discretization used to obtain the FOM by discretizing~\eqref{eq:strong_NSE} in space. To form the reduced basis, we first collect the FOM solution at $N_t\in\mathbb{N}$ time instances, or snapshots, over a training time horizon $[t_i,t_t]$, where $t_t<t_f$, and center the snapshots about their temporal mean $\umean(\bx)=\tfrac{1}{N_t}\sum_{j=1}^{N_t}\bu^h(t^{(j)},\bx)$. 
We then arrange the mean-centered fluctuations into the snapshot matrix $U\in\mathbb{R}^{N\times N_t}$ and compute the rank-$r$ POD basis $\PODbasis_r$. 
When the data are generated by an FE simulation, this discrete construction is equivalent to the corresponding continuous ROM projection~\cite{Gunzburger:2007,Kalashnikova:2014,volkwein2013proper}. 
Therefore, both representations are used interchangeably, writing $\basisfunc_i=\sum_{j=1}^N \phi^r_{ij}\boldsymbol{\varphi}^h_j\in\mathbf{X}^r$ for the ROM basis functions and $\basisvec_i$ for the columns of $\PODbasis_r$. 
This dual notation allows us to transfer the ROM differential filter introduced in Section~\ref{sec:diff_filtering}, which is naturally posed in the continuous setting~\cite{Moore2025_ADLROM,SANFILIPPO2023_ADL,Wells2017EF,xie2018-NA_Leray}, to the 
discretely formulated OpInf setting. 

Since POD modes are global, an arbitrary linear combination of them cannot match nonzero boundary values, and enforcing inhomogeneous Dirichlet conditions in a standard POD/Galerkin expansion is therefore nontrivial. 
Centering addresses this issue: for time-independent Dirichlet data, all snapshots satisfy the same boundary conditions, and therefore their mean $\umean(\bx)$ does as well. Therefore, the fluctuations $\bu^h(t_i) - \umean(\bx)$, together with every mode built
from them, vanish on the boundary. %
The ROM~\eqref{eq:state_var} then satisfies the boundary conditions strongly, for every $t$ and every fluctuation %
\cite{Sw19}. 
This is boundary enforcement method~\uppercase\expandafter{\romannumeral 1\relax} of~\cite{Gunzburger:2007}: the centering term carries the non-homogeneous boundary values, while the modal correction vanishes on the boundary.  This leads to the ROM approximation 
\begin{equation}
\romfunc(t, \bx) = \umean(\bx) + \sum_{i=1}^r \statecoef_i(t) \basisfunc_i(\mathbf{x}),
\qquad \state(t) = [\statecoef_1(t), \dots, \statecoef_r(t)]^T .
\label{eq:state_var}
\end{equation}
Since $\umean$ and each $\basisfunc_i$ are linear combinations of discretely divergence-free FOM data, incompressibility holds discretely within the trial space and the pressure term drops out~\cite{rozza2022advanced}. 
Substituting $\romfunc$ into~\eqref{eq:strong_NSE} and testing against $\basisfunc_i$ gives the intrusive quadratic G-ROM,
\begin{equation}
    M \frac{d \state }{d t} = \GROMconstantOP + \GROMlinearOP \state +  \state^T \GROMquadOP \state,
    \label{eq:GROM_dynamical_system}
\end{equation}
with operators
\begin{align*}
M_{ij} &= (\basisfunc_i, \basisfunc_j), \qquad
    \GROMconstantOP_i = (\basisfunc_i, f) - (\basisfunc_i, \umean \cdot \nabla \umean) - \mathrm{Re}^{-1} (\nabla \basisfunc_i, \nabla \umean), \\
    \GROMlinearOP_{ij} &= -(\basisfunc_i, \umean \cdot \nabla \basisfunc_j) - (\basisfunc_i, \basisfunc_j \cdot \nabla \umean) - \mathrm{Re}^{-1} (\nabla \basisfunc_i, \nabla \basisfunc_j), \qquad
    \GROMquadOP_{ijk} = -(\basisfunc_i, \basisfunc_j \cdot \nabla \basisfunc_k).
\end{align*}
OpInf postulates an equivalent quadratic structure,
\begin{equation}
\frac{d \state }{d t} = \constantOP + \linearOP \state + \quadOP ( \state \otimes \state),
\label{eq:Opinf_dynamical_system}
\end{equation}
but infers the reduced operators $\constantOP\in\mathbb{R}^{r}$, $\linearOP\in\mathbb{R}^{r\times r}$, $\quadOP\in\mathbb{R}^{r\times r^2}$ from data rather than from the FE operators.

\begin{remark}
\label{rem:mass_matrix}
Note that no OpInf analogue of the ROM mass matrix $M$ in~\eqref{eq:GROM_dynamical_system} appears in~\eqref{eq:Opinf_dynamical_system}. 
This is consistent with standard projection-based ROMs when the basis $\PODbasis_r$ is orthonormal with respect to the FE mass matrix and, consequently, the $L^2$ inner product in $\mathbf{X}^h$. 
When the FE mass matrix is unavailable, however, such a basis cannot be constructed, and the Euclidean POD basis will not in general have this property. 
In that case, multiplying~\eqref{eq:GROM_dynamical_system} by $M^{-1}$ yields a system of the form~\eqref{eq:Opinf_dynamical_system} with operators $M^{-1}\GROMconstantOP$, $M^{-1}\GROMlinearOP$, and $M^{-1}\GROMquadOP$. 
The two systems have the same trajectories, but they pose different regularized least-squares problems (with different symmetries) and therefore yield different inferred operators. 
In this paper, we use the completely non-intrusive POD construction described above for the CDR and FPC test cases, which corresponds to a (formal) multiplication by $M^{-1}$. Conversely, we use a function space formulation \cite{QFW21} for the MCF test case, in which the inner product is directly specified.
\end{remark}

\subsubsection{The OpInf learning problem and its regularization}
    \label{sec:opinf_regularization}

Let $\widehat{U}=\PODbasis_r^T U \in \mathbb{R}^{N_t \times r}$ denote the reduced snapshot matrix, whose $i$-th column is $\state(t_i) \in \mathbb{R}^{r}$, and let $D_t(\widehat{U})$ contain the corresponding time derivatives,
estimated by, e.g., finite differences when unavailable from the FOM. 
The OpInf learning problem computes the reduced model operators via
\begin{equation} \label{eq:opinf_learning}
    \min_{\constantOP, \linearOP, \quadOP}
    \bigl\| D_t(\widehat{U}) - \constantOP \mathbf{1}^T - \linearOP \widehat{U} - \quadOP \widehat{U}^{*2}  \bigr\|_F^2
    + \beta_1 \bigl(\| \constantOP \|_F^2 + \| \linearOP  \|_F^2\bigr) + \beta_2 \| \quadOP \|_F^2,
\end{equation}
where $\mathbf{1} \in \mathbb{R}^{N_t}$ is the vector of ones, $\widehat{U}^{*2}$ is the column-wise Kronecker (i.e, Khatri-Rao) product of $\widehat{U}$, and $\hyperparameter=[\beta_1,\beta_2]^T$ are the $L^2$ or Tikhonov regularization hyperparameters. 
As written, the quadratic operator $\quadOP \in \mathbb{R}^{r \times r^2}$ in~\eqref{eq:Opinf_dynamical_system} is not uniquely determined: the symmetry $\statecoef_j \statecoef_k = \statecoef_k \statecoef_j$ makes the
columns of $\widehat{U}^{*2}$ redundant, so distinct operators produce identical dynamics
and the data matrix in~\eqref{eq:opinf_learning} is rank deficient. 
We therefore eliminate the redundant degrees of freedom and learn an operator of dimension $r \times r(r+1)/2$, as described in %
\cite{Kramer2024_Opinf}.

Following~\cite{McQuarrie2021Regularized,QFW21}, Tikhonov regularization is added to~\eqref{eq:opinf_learning} to reduce overfitting and to account for model misspecification or other sources of error, with scalar regularization hyperparameters.
The constant and linear terms are regularized separately from the quadratic term, which improves performance at the cost of a larger hyperparameter search: a grid with $n$ candidate values per parameter requires $n^2$ solves of~\eqref{eq:opinf_learning}.
For a broader perspective, we consider three strategies in the present paper:
\begin{enumerate}
\renewcommand{\theenumi}{\Roman{enumi}}
    \item \emph{Unregularized OpInf}: $\hyperparameter=\mathbf{0}$; viable for simple problems, but inadequate for the flows considered here. \label{model:unregularized}
    \item \emph{Weakly regularized OpInf}: $\hyperparameter$ minimizes the training error $\|\widehat{U} - \widetilde{U}^{\,\mathrm{train}}(\hyperparameter)\|_F$ over the search grid. This improves the performance of unregularized OpInf but generally leads to overfitting.
    \label{model:weakly_regularized}
    \item \emph{Strongly regularized OpInf}: $\hyperparameter$ minimizes the same training error among candidates that do not exhibit excessive coefficient growth over a time horizon at least as long as the prediction horizon, with $N_k$ time instants, based on the strategy proposed in %
    \cite{McQuarrie2021Regularized}. \label{model:strongly_regularized}
\end{enumerate}
For model~\ref{model:strongly_regularized}, we measure the growth of the ROM approximate solution, $\widetilde{U}(\hyperparameter)$, against the temporal mean
$\widehat{U}^{\mathrm{mean}}=\frac{1}{N_t+1}\sum_{i=0}^{N_t}\widehat{U}_i$ of the reduced training data, discarding every $\hyperparameter$ whose growth ratio
\begin{equation}
    G(\hyperparameter) = \frac{\max_{N_t< i \leq N_k} \| \widetilde{U}_i(\hyperparameter) - \widehat{U}^{\mathrm{mean}} \|_\infty}
    {\max_{0\leq i \leq N_t} \| \widehat{U}_i - \widehat{U}^{\mathrm{mean}} \|_\infty},
    \label{eq:opinf_growth_ratio}
\end{equation}
exceeds a prescribed tolerance $\tau$ (e.g., $\tau=1.2$) as likely to yield an unstable model.
This criterion uses only ROM solutions beyond the training horizon and requires no FOM data in that regime, so it remains computationally inexpensive. 
Once $\hyperparameter$ is selected through the grid search, which is %
parallelizable~\cite{FARCAS2025109619}, the learned operators are frozen and used throughout the entire online phase.
Further extensions to OpInf refine either the regularization, via Bayesian
formulations~\cite{GUO2022bayesianopinf,mcquarrie2025activelearningdatadrivenreduced}, Gershgorin-based eigenvalue shifts~\cite{GKIMISKIS2025localizedopinf}, or parametric extensions~\cite{farcas2023parametric,mcquarrie2023parametric}, or the structure of the learned operators, through Hamiltonian and energy-conserving formulations~\cite{Gruber2025Hamiltonian,SHARMA2022Hamiltonian} and nested constructions that reduce sensitivity to $\hyperparameter$~\cite{Aretz:2025}.

The learning problem~\eqref{eq:opinf_learning} minimizes a residual assembled instant-by-instant from the training data, with no coupling across time. The ROM is then deployed online by integrating the learned operators over the time domain of interest. 
Therefore, a small residual does not necessarily imply an accurate reduced solution, which suggests selecting $\hyperparameter$ based on the reduced solution $\widetilde{U}(\hyperparameter)$ of the OpInf ROM~\eqref{eq:Opinf_dynamical_system}.

\begin{remark}[Regularization for parametric settings]
    The regularization procedure described in this section is formulated primarily for predictive ROMs, namely OpInf models used to extrapolate beyond the time interval on which they are trained.
    This setting is employed predominantly in our numerical experiments and provides the main motivation for the proposed EFR-OpInf method.
    The same general principles can also be extended to parametric OpInf models, for which the regularization parameters must account for accuracy and stability across multiple parameter instances.
    In Section~\ref{sec:CDR}, we consider a predictive problem that additionally incorporates parametric variation.
    Further details on the regularization of parametric OpInf models can be found in
    \cite{farcas2023parametric,mcquarrie2023parametric}.
\end{remark}

Most of the OpInf literature 
partitions the time domain of interest into \emph{training} and \emph{prediction} intervals, noting that selection of the regularization hyperparameters may involve the imposition of 
bound constraint on the learned coefficients over a time horizon that extends beyond the final prediction time~\cite{McQuarrie2021Regularized}.
To more rigorously assess the robustness of the considered filtering strategies, we instead partition the time domain of interest into three subintervals: \emph{training}, \emph{validation}, and \emph{testing}.
We explain this choice in more details below.

\vspace{1pc}
{\noindent\bf Training, validation, and testing.}
The terms \emph{training}, \emph{validation}, and \emph{testing} are used here in the usual machine-learning sense.
For all standard OpInf models, 
only high-fidelity data from the training interval are supplied to the regression problem~\eqref{eq:opinf_learning}.
For weakly regularized OpInf, model~\ref{model:weakly_regularized}, the regularization hyperparameters $\hyperparameter$ are selected solely by minimizing the training error.
Consequently, neither the validation nor the test data influence the learned operators or the choice of $\hyperparameter$, and both intervals may be used to assess the out-of-sample accuracy and stability of the resulting ROM.
For strongly regularized OpInf, model~\ref{model:strongly_regularized}, the hyperparameters are selected by minimizing the training error subject to the boundedness constraint on the learned reduced coefficients described above, which is enforced over the validation interval.
Thus, the validation interval is used for hyperparameter selection, whereas the test interval is reserved exclusively for evaluating the accuracy and stability of the selected ROM on previously unseen data.
Much of the OpInf literature instead divides the time domain into only \emph{training} and \emph{prediction} intervals.
Within this framework, a training--validation--test split can be introduced by dividing the prediction interval into a validation region, used to select the regularization hyperparameters, and a subsequent test region, used for evaluation on unseen conditions.
One could alternatively enforce the boundedness constraint over the entire combined validation-and-test horizon.
However, for short training horizons as those considered in this work, enforcing boundedness over a long extrapolation interval may require excessively large regularization parameters.
This can drive the learned operators toward small values and produce overly dissipative or otherwise inaccurate predictions.
Conversely, enforcing boundedness only over the shorter validation interval may yield models that remain stable during validation but become unstable over the longer test horizon.
This tradeoff between short-term accuracy and long-term stability is one of the primary motivations for the filtering extension proposed in this paper.
The EFR-OpInf models are constructed from standard OpInf ROMs.
A limited amount of high-fidelity data from the validation interval are used to select the filtering parameters and assess stability over that interval.
The resulting EFR-OpInf models are evaluated over the previously unseen test interval to determine whether ROM solution stability persists over longer prediction horizons.
For all three
test problems, the coefficient boundedness constraint is enforced over the
validation interval to select the filtering parameters. This constraint
requires only the predicted ROM coefficients and does not use high-fidelity
validation data. For the MCF problem, one representative FOM solution from
the validation interval is additionally used to identify an initial set of
candidate filtering parameters. This preliminary step is not essential, as
the filtering parameters could instead be selected solely by enforcing the
coefficient boundedness constraint. No high-fidelity data from the test
interval is used for model or parameter selection. The resulting EFR-OpInf
models are then evaluated over the previously unseen test interval to
determine whether their accuracy and stability persist over longer prediction
horizons. %

\subsection{ROM Filtering Algorithms}
    \label{subsec:filtering_methods}

Spatial filtering is an LES strategy that has traditionally been applied at a FOM level to overcome issues of limited mesh resolution and complex turbulent interactions~\cite{berselli2006LES,Sagaut2006LES}. 
Spatial filtering has also been applied at the ROM level to mitigate spatial under-resolution, which can lead to a non-physical accumulation of energy in the resolved ROM modes~\cite{Wells2017EF}.
These problems seen in G-ROMs can also be observed in OpInf ROMs with limited training data.
Indeed, as the numerical results of Section \ref{sec:numerical_results} will show, for short training horizons, choosing OpInf hyperparameters that ensure the %
model's stability beyond the training interval can yield inaccurate predictions. 
Conversely, choosing OpInf hyperparameters that strictly minimize the training error can yield unstable models for long-term prediction.
To address these challenges, we propose a novel strategy: the incorporation of online, ROM-level spatial filtering into the OpInf approximate solution to quell spurious numerical oscillations and promote stability in scenarios with limited amounts of training data. 
As we will show in Section~\ref{sec:numerical_results}, the proposed ROM spatial filtering strategy reduces, and in some cases even eliminates, the need for standard regularization in the OpInf process, while introducing stabilization parameters with a clearer physical interpretation.  

The proposed EFR-OpInf approach presented in Section \ref{sec:EFR_OpInf} will use the following three ROM spatial filters:
(i) the ROM projection filter (Section \ref{sec:Proj_Filter});
(ii) the ROM differential filter (Section \ref{sec:diff_filtering}), and 
(iii) the ROM hybrid projection-differential filter (Section \ref{sec:hybrid_proj_diff_filter}). 
We note that the first two filters have been extensively used in the development of ROM closures and stabilizations, %
whereas the third filter was recently introduced in %
\cite{strazzullo2025variational}.
Moreover, because each filter has proved effective in different modeling and computational settings, all three are viable candidates for EFR-OpInf.
We therefore assess the performance of all three filters numerically in Section~\ref{sec:numerical_results}.  To %
our knowledge, no comprehensive comparison of these ROM filters is currently available. 

These three ROM filters are presented in the context of fluid flow, where the velocity field is separated into mean and fluctuating components. 
Following~\cite{Wells2017EF}, all ROM spatial filters are applied only to the small-scale fluctuating components.  
Thus, the filtered ROM velocity $\filfunc$, with separated mean and fluctuating components, can be written as follows: 
\begin{align}
    \filfunc(t,\bx) &= \umean(\bx) + \sum_{i=1}^r \filcoef_i(t) \basisfunc_i(x), \;\;\;\;\; \filstate = [\filcoef_1(t), \dots, \filcoef_r(t)]^T. \label{eq:filter_var}
\end{align}
Observe that both the ROM velocity field in \eqref{eq:state_var} and the filtered ROM velocity field in \eqref{eq:filter_var} are written in terms of the same basis functions, but potentially with different coefficients. 

\subsubsection{ROM Projection Filter}
    \label{sec:Proj_Filter}

The ROM projection filter is a simple, fully discrete, and non-intrusive spatial filter.
It has been successfully used in the development of ROM closures and stabilization strategies for convection-dominated, including turbulent, flows~\cite{kaneko2022LerayROM,KANEKO2020POD_ROM, Koc2022Verifiability, Reyes2025Verifiability, WANG2012PODClosure,Wells2017EF}.
The filter maps a ROM state represented in an $r$-dimensional reduced space onto a coarser $r'$-dimensional space, with $r' < r$, by retaining only the first $r'$ vectors of the POD basis $\PODbasis_r$.
This construction is motivated by the energy ordering of the POD modes: in many convection-dominated flow problems, the leading modes predominantly represent the large-scale, energetic structures, whereas higher-index modes tend to contain less energetic and often smaller-scale features~\cite{holmes2012turbulence}.
Consequently, projecting the ROM state onto the space spanned by the leading $r'$ POD modes removes part of the small-scale content and provides a natural ROM-level spatial filtering mechanism. %

Next, we explain the construction of the ROM projection.
Given an $r$-dimensional POD basis $\PODbasis_r$, let $\PODbasis_{r'}$ be the restriction of $\PODbasis_{r}$ to its first $r'<r$ columns. 
Recalling that the projection of the mean-centered snapshot matrix into the $r$-dimensional space determined by the basis is given by $\PODbasis_r \PODbasis_r^T U$, further projection into the $r'$-dimensional space is given by $\PODbasis_{r'} \PODbasis_{r'}^T \PODbasis_r \PODbasis_r^T U$. 
Note that, due to the orthonormality of $\PODbasis_r$, the column-wise inner products satisfy
$
\PODbasis_{r'}^T \PODbasis_r =
\begin{bmatrix}
I_{r'} &  \mathbf{0}_{(r') \times (r-r')} 
\end{bmatrix}.
$
This implies that $\PODbasis_{r'}\PODbasis_{r'}^T\PODbasis_r\PODbasis_r^T U$
performs the following sequence of operations:
$\PODbasis_r^T U$ computes the $r$-dimensional reduced coefficients of $U$;
$\PODbasis_{r'}^T\PODbasis_r(\PODbasis_r^T U)$ retains only the first $r'$
components of these coefficients; and 
left multiplication by
$\PODbasis_{r'}$ lifts the truncated coefficient vector back to
$\mathbb{R}^N$ as a linear combination of the first $r'$ basis vectors.
Thus, projection onto the $r'$-dimensional reduced space is equivalent to
truncating the final $r-r'$ reduced coefficients.

The ROM projection filter is typically performed at a ROM level. When applied to an 
unfiltered set of coefficients $\state$, the filtered coefficients $\filstate$ may be computed as,
\begin{equation}
    \begin{bmatrix}
        u_1, &
        u_2, &
        \dots, & &
        u_r
    \end{bmatrix}^T = \state \mapsto \filstate =     \begin{bmatrix}
        u_1, &
        \dots &
        u_{r'}, &
        0, &
        \dots, &
        0
    \end{bmatrix}^T. 
    \label{eq:ROM_Proj_Filter}
\end{equation}
It is generally more convenient in software implementations to consider the projection filter as setting the excess $r - r'$ coefficients to 0, as in~\eqref{eq:ROM_Proj_Filter}, rather than eliminating the excess basis functions.  This way, both $\state$ and $\filstate$ live in $\mathbb{R}^n$ and are compatible with the same vector operations. 
We note that both the coefficient and basis perspectives of this projection filter produce equivalent numerical results. 

\begin{remark}
    The ROM projection filter can also be formulated in terms of a continuous projection analogous to the continuous POD~\cite{volkwein2013proper}.
    We have chosen to demonstrate and implement the discrete version of the projection filter
    since it can be easily implemented into existing OpInf %
    software without any alteration or additional data assumptions. 
\end{remark}

\subsubsection{ROM Differential Filter}
    \label{sec:diff_filtering}

Just as the ROM projection, the ROM differential filter has been used to construct ROM closures and stabilizations~\cite{strazzullo2021_EFRConsistency,Wells2017EF,xie2018-NA_Leray}.
The ROM differential filter is applied to the fluctuating components of ROM functions $\romfunc$ given in equations \eqref{eq:state_var}.  %
Specifically, for any known ROM velocity $\romfunc$, the coefficients $\filcoef$ of the filtered velocity $\filfunc$ in \eqref{eq:filter_var} can be obtained by %
solving the system 
\begin{equation*}
    ( I - \delta^2 \Delta)  \left( \sum_{i=1}^r \filcoef_i(t) \basisfunc_i(x) \right) =    \sum_{i=1}^r \statecoef_i(t) \basisfunc_i(x).
\end{equation*}
Testing against the ROM basis $\basisfunc_j$, for $j = 1,\dots r,$ and noting that the basis functions $\basisfunc$ have been specifically constructed with homogeneous Dirichlet boundary conditions, this is equivalent to the weak-form system 
\begin{equation}
   \sum_{i=1}^r \filcoef_i(t)\left(  \basisfunc_i, \basisfunc_j \right) + \delta^2 \sum_{i=1}^r \filcoef_i(t) \left( \nabla \basisfunc_i, \nabla \basisfunc_j \right) =  \sum_{i=1}^r \statecoef_i(t) \left( \basisfunc_i, \basisfunc_j \right), \; j = 1, \dots, r.  
   \label{eqn:rom-df}
\end{equation}
Recalling that the ROM mass matrix, $M$, has components $M_{ij} = (\basisfunc_i, \basisfunc_j)$, and defining the stiffness matrix $S$ to have components $S_{ij} = (\nabla \basisfunc_i, \nabla \basisfunc_j)$, 
the ROM vector of filtered coefficients $\filstate$ is compactly expressed as the solution to the linear system 
\begin{equation}
    (M + \delta^2 S) \filstate = M \state. 
    \label{eq:ROM_DF_Filter}
\end{equation}

The ROM differential filter is attractive in traditional %
G-ROM settings partly because it can be implemented entirely in reduced coefficients without requiring explicit FOM field reconstructions, in contrast to Gaussian or other convolution-based filters. %
Additionally, the effects of the ROM differential filter are not as closely linked to the POD basis as in the ROM projection filer,
allowing for finer control over the filtering effect compared to 
the ROM projection filter. This flexibility comes at the cost of requiring the ROM-level differential operators $M$ and $S$, and therefore makes the approach somewhat more involved than a fully data-driven OpInf workflow. The $M$ and $S$ operators may be obtained directly from the FOM code, when such access is available, or constructed independently from the full-order mesh. The latter does not require access to the FOM implementation, but does require additional information beyond the snapshot data.
The inclusion of ROM differential filtering in this work 
helps evaluate the utility of standard ROM filtering approaches in the OpInf setting. We aim to quantify whether the positive features of the filter justify the additional access to FOM-level information required to assemble \(M\) and \(S\), and to assess whether developing fully non-intrusive alternatives with properties analogous to the ROM differential filter is beneficial for OpInf ROMs. 

\subsubsection{%
Hybrid ROM Projection-Differential Filter} %
    \label{sec:hybrid_proj_diff_filter}

While many LES models are based on differential filtering at the FOM level~\cite{berselli2006LES,Germano1986Differential, Sagaut2006LES} and at the ROM level~\cite{GIRFOGLIO2023FILTER, strazzullo2021_EFRConsistency,Wells2017EF, xie2018-NA_Leray}, a common %
challenge for these models is their potential to be overly diffusive. This has led to the development of less diffusive approximate deconvolution FOM~\cite{Layton2012AD} and ROM~\cite{Moore2025_ADLROM} models. 
In this paper, we address this 
concern with a variational multiscale filter recently developed in~\cite{strazzullo2025variational}. %
Denoted as the hybrid ROM projection-differential filter, this procedure separates the ROM basis basis functions into large, dominant structures and smaller, oscillatory structures, akin to the projection filter. That is, defining some $r' < r$, we assume basis functions $1,\dots r'$ represent fully resolved components, while basis functions $r'+1, \dots, r$ represent marginally resolved components. The hybrid filter %
keeps the first $r'$ modes unchanged and applies differential filtering to modes $r' +1$ through $r$. The key idea here is to leave the dominant and typically smooth structures (that do not require filtering) intact while extracting as much information as possible out of the small-scale components. 
For clarity, the precise definition of the hybrid projection-differential filter used in Section~\ref{sec:numerical_results} is provided below. %

\begin{definition}[Hybrid ROM Projection-Differential Filter]
\label{def:HPDF}
Let $\basisfunc_i(\bx)$, $i = 1, \dots, r$, be a basis for $\mathbf{X}^r$, and let $F_{DF}(\cdot)$ refer to the output of the differential filter given by~\eqref{eqn:rom-df}. %
Given any ROM function $\romfunc(t, \bx)$, the hybrid projection differential filter $F_{HDF}(\romfunc(t, \bx))$ is:
\begin{equation}
F_{HDF}(\romfunc(t, \bx)) %
= \umean(\bx) + \sum_{i=1}^{r'} \statecoef_i(t) \basisfunc_i(x) + F_{DF}\left(\sum_{i=r'+1}^r \statecoef_i(t) \basisfunc_i(x)\right).
\label{eq:partial_differential_filtering}
\end{equation}
\end{definition}

The cut-off dimension $r'$ used in the above definition is similar to that of the projection filter. 
However, here it does not indicate truncation, instead separating the 
filtered components from the unfiltered components of the ROM solution.
Note that 
any ROM-level filter could be used in place of $F_{DF}(\cdot)$ in~\eqref{eq:partial_differential_filtering}. We have chosen 
to use the differential filter due to its success and widespread use within the ROM literature. Since the ROM differential filter has been noted to be over-diffusive~\cite{Moore2025_ADLROM,SANFILIPPO2023_ADL}, we believe this hybrid filter may help address this point. %

In terms of vector representations, if $\state$ is the vector of ROM coefficients and $F_{HDF}(\cdot)$ is the hybrid projection differential filter described in Definition \ref{def:HPDF}, 
then its action on coefficients is given by
\[
F_{HDF}(\state) = [\statecoef_1, \dots, \statecoef_{r'}, \filcoef_{r'+1}, \dots \filcoef_{r}]^T,
\]
where ${\filcoef_{i}}$ denotes the $i^{\text{th}}$ component of the output of the fully discretized ROM differential filter~\eqref{eq:ROM_DF_Filter}. 
That is, we obtain the last $r - r'$ coefficients of the hybrid filter as the last $r - r'$ coefficients of the $r-$dimensional differential filter~\eqref{eq:ROM_DF_Filter}, while leaving the first $r'$ components untouched. We briefly note that this procedure could be made more efficient by solving only the $(r-r')$-dimensional differential filter, rather than the complete $r$-dimensional differential filter, but this cost is typically negligible in practice.

\medskip

We conclude this section by reiterating that the projection filter requires no additional information beyond what is needed for the standard POD-based OpInf model.
Both the differential filter and the new hybrid projection-differential filter require the solution of an $r$-dimensional linear system. 
However, this generally incurs negligible computational cost, since $r$ is $\mathcal{O}(10)$ for the ROM problems considered here, even when the underlying FOM has thousands to millions of DoFs.

\section{Evolve-Filter-Relax OpInf (EFR-OpInf)}
\label{sec:EFR_OpInf}

To reduce the sensitivity of OpInf models to problem-dependent regularization parameters and alleviate the resulting tradeoff between short-term accuracy and long-term stability, we introduce a filtering-based OpInf framework that supplements the learned reduced dynamics with ROM-level spatial filtering. 
Rather than relying exclusively on Tikhonov regularization to simultaneously control model accuracy and long-term stability, the proposed approach uses filtering to suppress unresolved oscillations and stabilize the online evolution. 
This separation provides a more flexible mechanism for balancing accuracy and stability, and can reduce the need for large regularization parameters that can produce overly dissipative or otherwise inaccurate predictions. %

Evolve-Filter-Relax is a simple and modular LES-inspired model that has frequently been used to stabilize fluid simulation at both a FOM level~\cite{Bertagna2016_EFR,Mullen1999Filtering, strazzullo2021_EFRConsistency} as well as a ROM level~\cite{ GIRFOGLIO2023FILTER,IVAGNES2026_EFR, strazzullo2021_EFRConsistency, TSAI2025_timerelaxation, Wells2017EF}. 
The key idea of EFR is to moderate the growth of solutions by incorporating just enough spatial filtering to modulate spurious oscillations caused by under-resolution and instabilities in the time-integration process. 
To perform the three-step process of EFR, we require a time integration scheme (e.g., forward Euler) applied to the OpInf dynamical system~\eqref{eq:Opinf_dynamical_system}, along with a spatial filtering method for the evolved coefficients, e.g., the ROM projection filter~\eqref{eq:ROM_Proj_Filter} or ROM differential filter~\eqref{eq:ROM_DF_Filter}. The EFR algorithm then proceeds according to the online portion of Algorithm \ref{alg:efr-opinf}.

\begin{algorithm}
\caption{EFR-OpInf}
\label{alg:efr-opinf}
\begin{algorithmic}[1]

\Statex \vspace{0.5\baselineskip}
\State \textbf{OpInf Construction (Offline)}
\Require Training window snapshot data \(U \in \mathbb{R}^{N\times N_t}\), target energy \(\eta\) (e.g., \(0.999\)), time-integration scheme $E$ (e.g., RK4), set of potential hyperparameter values $\hyperparameter_i$ for $i = 1, \dots, N_p$ with  $N_p \in \mathbb{N}$, (Optional) maximum growth factor $M$ (see equation \eqref{eq:opinf_growth_ratio}).
\State Extract POD basis \(\PODbasis\) from snapshots \(X\); set reduced dimension \(r\) according to target energy \(\eta\).
\State \(\PODbasis_r \leftarrow \PODbasis(:,1\!:\!r)\)
\State Compress snapshots \(U \mapsto \widehat{U}\) using $\PODbasis_r.$
\State $\texttt{current error} \gets \infty$
\For{$i = 1, \dots, N_p$}
\State Solve OpInf learning problem~\eqref{eq:opinf_learning} for parameter value $\hyperparameter_i$. 
\State Perform online run of OpInf using time-integration scheme $E$ to produce OpInf solution $\widetilde{U}(\hyperparameter_i)$.
\State (Optional) Compute coefficient growth factor $G$, disqualify $\hyperparameter_i$ if $G > M$ (equation~\eqref{eq:opinf_growth_ratio}).
\State Calculate OpInf training error $e({\hyperparameter_{i})} = \|\widehat{U} - \widetilde{U}^{train}(\hyperparameter_i)\|_F$
\If{$e({\hyperparameter_{i}}) < \texttt{current error}$}
\State \(\hyperparameter_{opt} \gets \hyperparameter_i\)  
\State $\texttt{current error} \gets e({\hyperparameter_{i}})$
\EndIf
\EndFor
\State \Return OpInf operators parameterized by $\hyperparameter_{opt}$ and resulting time-evolution scheme \(E(t, \state^{n}; \hyperparameter_{opt}) \mapsto \hat{u}^{n+1}\).

\Statex \vspace{0.5\baselineskip}
\State \textbf{EFR-OpInf (Online)}
\Require Prediction time span \((t_0,t_f)\),
timestep \(\Delta t\), OpInf time-evolution scheme \(E(t, \state^{n}; \hyperparameter) \mapsto \hat{u}^{n+1}\), filter method \(F(\state) \mapsto \filstate\), relaxation parameter $\chi \in (0,1]$.
\State $n \gets 0$
\State $t \gets t_0$

\While{\(t < t_f\)}
  \State \textbf{Evolve}: \(\hat{w}^{n+1} \leftarrow E(t, \state^{n}, \hyperparameter)\)
  \State \textbf{Filter}: \(\overline{\hat{w}}^{\,n+1}\leftarrow F(\hat{w}^{n+1})\)
  \State \textbf{Relax}: \(\state^{n+1} \leftarrow \chi\,\overline{\hat{w}}^{\,n+1} + (1-\chi)\,\hat{w}^{n+1}\)
  \State $t \gets t + \Delta t$
  \State $n \gets n + 1$
\EndWhile
\end{algorithmic}
\end{algorithm}

The relaxation parameter, $\chi \in [0,1]$, demarcates the closely related Evolve-Filter methods (e.g., \cite{Wells2017EF}) from Evolve-Filter-Relax methods (e.g., \cite{gunzburger2019EFR_ROM, strazzullo2021_EFRConsistency}), and has been shown~\cite{Bertagna2016_EFR} to be important in limiting the amount of diffusion applied to the solution. 
The convex combination, parameterized by $\chi$, determines the amount of filtering included in the ultimate solution update, with larger values of $\chi$ corresponding to greater amounts of filtering.
Numerical analysis at a FOM level~\cite{neda2012_EFR_analysis} suggests the choice $\chi \approxeq \mathcal{O}(\Delta t)$ to limit added dissipation. 

As mentioned, we propose combining the OpInf framework with the EFR strategy in order to increase OpInf's stability and accuracy, reducing the need for large regularization parameters that yield inaccurate results.
The resulting model is denoted \emph{evolve-filter-relax OpInf (EFR-OpInf)} and summarized in Algorithm \ref{alg:efr-opinf}.
Note that EFR is modular, making it highly complementary to OpInf. 
 Furthermore, while the evolve step in EFR-OpInf inherently relies on the purely data-driven OpInf model, the filtering step can be adapted to the FOM information available in a given application. The ROM projection filter requires only the snapshot data already needed for OpInf, whereas the ROM differential filter additionally requires the assembly of simple, standard operators commonly used in FE simulations. When further FOM information or capabilities are available, more recent ROM filtering approaches may also be incorporated, such as high-order algebraic filters~\cite{TSAI2025_timerelaxation}. Alternatively, fully data-driven filters may be employed, such as the recently introduced StabOp ROM~\cite{tsai2026stabopdatadrivenstabilizationoperator}.

Filtering methods have not yet been widely studied for OpInf models. To our knowledge, the only approach that employs filtering in an OpInf setting is found in~\cite{farcas2023filtering}, and uses a FOM-level spatial filter to smooth the training data before training an OpInf model on this filtered data. This offline filtering approach is more similar to the OpInf regularization methods discussed in Section \ref{sec:opinf_regularization} than to the online EFR-OpInf model that we propose here, since its motivation is to improve
the training of OpInf models. The novelty of our approach is twofold. First, we employ ROM-level and not FOM-level filtering. Second, 
the use of ROM-level information
enables the use of filtering in the online stage of the ROM
as opposed to during offline training. Instead of improving the training procedure of an OpInf ROM, our method improves the performance of 
one that is already deployed.

The main advantage of the proposed EFR-based online filtering approach, summarized in Algorithm \ref{alg:efr-opinf}, is that this filtering can be turned on and off and (in most cases) adjusted at will. While a standard pre-trained OpInf model carries a certain
amount of regularization 
baked into the learned operators that remains fixed for the entire simulation, 
online EFR parameters such as $\chi$ can be adjusted at will in the online stage. 
The time-predictive EFR-OpInf in Section \ref{sec:numerical_results} %
employs filtering only 
outside the training regime, thus 
obtaining %
accuracy within the training regime while remaining stable in the predictive regime. 
We note that filters of varying strength and intrusiveness can be selected based on 
information available in the problem setting at hand, and are 
not restricted to those presented in %
Section~\ref{subsec:filtering_methods}. 

A second important advantage of EFR-OpInf comes from the \emph{physical interpretability} of the EFR filtering approach. In standard OpInf, just as in machine learning, the regularization parameter is a property of the operators themselves, 
and it is difficult to 
quantify the 
effect of one set of parameters $\hyperparameter$ 
as opposed to another. 
For example, one could ask the following questions. 
What is the difference between choosing a parameter $\beta_1=10^{-5}$ instead of $\beta_1=10^{7}$ in OpInf? %
How would this choice be reflected in the resulting OpInf dynamical system?
Would the OpInf solutions be smoother or more oscillatory?
We emphasize that, in contrast, each of the EFR filters under consideration can be clearly visualized and has physical meaning built into its construction (see Appendix~\ref{appendix:spatial-filtering}), making them more interpretable. 
For example, the cutoff dimension for the ROM projection filter, $r'$, corresponds directly to retained energy in the ROM modes. Similarly, the ROM differential filter parameter, $\delta$, corresponds to a lengthscale that has clear physical meaning at a FOM level~\cite{berselli2006LES} and had recently been the subject of investigation at a ROM level~\cite{MOU2023_lengthscale}. 
We also note that there is rigorous mathematical support available for these filters.
For example, the differential filter 
was %
proved to be stable~\cite{IVAGNES2026_EFR,xie2018-NA_Leray}.  Moreover, parameter scalings based on the lengthscale $\delta$ were derived (and computed) for a Leray-type ROM using a differential filter in~\cite{xie2018-NA_Leray} and for a time relaxation ROM in~\cite{reyes2025timerelaxationNA}. %
Additional proofs on the stability of each of these filter methods %
were provided in \cite{moore2026filtering}. 
These numerical analysis results for the ROM filters will serve as a mathematical guide in choosing the EFR-OpInf parameters in the numerical investigation in Section~\ref{sec:numerical_results}.
Specifically, instead of searching the regularization parameters over several orders of magnitude (as is often done for the standard OpInf), we will search the filtering parameters in EFR-OpInf over only a few (i.e., $\mathcal{O}(10)$) values that are suggested by numerical analysis results.  

\section{Numerical Results}
\label{sec:numerical_results}

The proposed EFR-OpInf strategy is demonstrated in three scenarios: the 2D unsteady convection-diffusion-reaction equation with parametric dependence introduced through a multiplicative factor involving the inverse P\'eclet number, the two-dimensional transitional flow past a cylinder (FPC) at $\mathrm{Re}=100$, and the three-dimensional minimal channel flow at $\mathrm{Re} = $ 5,000.
Our numerical investigation will show that the under-resolved simulations of these convection-dominated (the first and third tests) and transitional (the second test)
flows can be challenging for the standard OpInf framework in the limited training data regime.
For each problem, we will compare EFR-OpInf models built using each of the proposed filtering methods: the ROM projection filter, the ROM differential filter, and the hybrid ROM projection-differential filter. 
These will be labeled as EFR-OpInf-Proj, EFR-OpInf-DF, and EFR-OpInf-HDF respectively.

Our filtered models will be compared against standard OpInf models using the regularization techniques~\ref{model:unregularized}--\ref{model:strongly_regularized} presented in Section \ref{sec:opinf_regularization}. 
In particular, the CDR case will consider unregularized as well as strongly regularized OpInf, which obtains the hyperparameter $\hyperparameter$ minimizing the training error and requires the coefficients of the ROM model to remain stable over a validation window.
For the FPC and MCF test cases, we consider the same strongly
regularized OpInf as well as a weakly regularized OpInf, in which 
the hyperparameters $\hyperparameter$ are selected solely by minimizing the training error. We do not consider
unregularized OpInf for the FPC and MCF problems because it produces unstable
solutions even within the training interval.

Because EFR stabilization is applied online to an already trained OpInf ROM, the performance of an EFR-OpInf model may depend on the regularization hyperparameter $\hyperparameter$ used to learn the underlying OpInf operators. 
To reduce this impact, for the CDR problem, ERF filtering is applied to 
the unregularized OpInf model (model~\ref{model:unregularized}); 
in this case, it is shown that EFR can eliminate the need for standard OpInf regularization. 
For the NSE problems, EFR filtering is applied to the weakly regularized OpInf model (model~\ref{model:weakly_regularized}). 
The strongly regularized OpInf model (model~\ref{model:strongly_regularized}) is included 
as a standard-practice baseline for OpInf approaches that do not employ online spatial filtering.
In our experience, weakly regularized OpInf models tend to overfit the training data: they achieve high in-sample accuracy but exhibit reduced stability and poorer generalization beyond their training interval. 
The EFR formulation can therefore leverage their accurate representation of the training dynamics while restoring stability in the extrapolation regime. 

We consider as metric the 
relative $L^2$ error at each timestep $t$, 
\[e_{L^2}(t) = \frac{\| \bu^h(t, \mathbf{x}) - \bu^r(t, \mathbf{x})\|_{L^2}}{\|\bu^h(t, \mathbf{x}) \|_{L^2}}. \] For the more complicated NSE problem setting, we will also compare the per-timestep kinetic energies of the ROM and FE velocity fields: defining $KE(\romfunc) = \frac{1}{2}\| \romfunc \|_{L^2}^2$, the relative kinetic energy error at time $t$ is given as,
\[e_{KE}(t) = \frac{| KE(\bu^h(t, \mathbf{x})) - KE(\bu^r(t, \mathbf{x}))|}{| KE(\bu^h(t, \mathbf{x})) |}. \]
For the %
FPC and MFC problems, we present an additional parameter and timing study for the various EFR-OpInf models in Section \ref{sec:NSE_FPC_Compare} (Figure \ref{fig:Pareto_Vary_Chi}) for FPC problem and Section \ref{sec:mfu_compare} (Figure \ref{fig:Pareto_Vary_Chi_MFU}) for the MCF problem. 
Both the CDR and %
the FPC simulations use FE data generated using the Python-based finite element platform FEniCSx~\cite{BarattaEtal2023,BasixJoss,ScroggsEtal2022}. 
The three-dimensional MCF simulations use SEM data generated with nekRS~\cite{fischer2022nekrs} and processed using
NekROM~\cite{kaneko_nekrom} to construct the POD basis.
The OpInf learning problem~\eqref{eq:opinf_learning} is solved using the Python package {\tt OpInf} \cite{opinf_python}. Finally, for models involving the differential filter, the requisite $M$ and $S$ from equation~\eqref{eq:ROM_DF_Filter} are computed by first obtaining the FOM mass and stiffness matrices from computation on the FOM mesh and then projecting these FOM operators onto the OpInf POD basis. %

\subsection{A High P\'eclet Number Parametric Convection-Diffusion-Reaction Problem}
    \label{sec:numerical-resutls-cdr}

The goal of this example is to demonstrate that EFR-OpInf can replace the standard regularized OpInf formulation with a fully interpretable spatial filtering strategy.
This is particularly useful for practitioners seeking a stable and computationally efficient model with minimal parameter tuning.
We consider a linear scenario containing distinct convection, diffusion, and reaction terms, while approximating its global dynamics with a linear OpInf model.
Given snapshots generated at a prescribed set of training parameter instances, the predictive task is to simulate the system trajectory at an unseen parameter configuration.

EFR-OpInf is used to stabilize the solution of the unregularized OpInf ROM (i.e., model~\ref{model:unregularized}).
For each of the three filters considered, we employ filter and relaxation parameter values motivated by numerical analysis, without performing an additional parameter search.
Since this example already incorporates parametric generalization, we evaluate future-state prediction only over the validation horizon, which is also used to select the regularization parameters for the strongly regularized OpInf model.
Following this procedure, all three EFR-OpInf variants remain stable and accurate over the entire prediction horizon without tuning the OpInf Tikhonov regularization parameters or conducting a separate search for the EFR parameters. 
The resulting workflow reduces implementation complexity and offline computational cost while providing a physically interpretable stabilization mechanism.

\label{sec:CDR}
\subsubsection{%
CDR FOM Setup}

The governing PDE for this example is given by the convection-diffusion-reaction equation
\begin{equation}
\begin{split}
    \ddt{u} - \epsilon \Delta u + \mathbf{b} \cdot \nabla u + \sigma u&= f, \; \text{in } \Omega \times (0, t_f],\\
    u|_{\partial \Omega} &= 0. 
    \label{eq:CDR_FOM}
\end{split}
\end{equation}
Here, $\Omega = [0,1] \times [0,1]$ is the unit square, $\mathbf{b} = [1,0]^T$, $\sigma = 1$, $f \equiv 1$, and the diffusion parameter $\epsilon$ is defined 
as the multiplicative inverse of the P\'eclet number, $\mathrm{Pe}$. 
Starting from an initial condition $u_0 \equiv 0$, the solution reaches a steady state after $t = 1.0$ (see Figure~\ref{fig:CDR_FOM_Steady}), forming a sharp boundary layer along the line $y = 1$. 
Accordingly, the FOM mesh is an evenly spaced grid with $h = 1/50 = 0.02$ for the majority of the domain, but is refined near the boundary to improve spatial resolution near the sharp boundary layer. 
The choice of quadratic elements leads to $11,000$ DoFs for the FE problem and produces a model able to resolve shocks with P\'eclet numbers up to $\mathrm{Pe}=10,000$. 
The FE problem is solved using the backward Euler method with timestep $\Delta t = 0.01$. 

\begin{figure}[ht]
    \centering
    \begin{subfigure}[t]{0.49\linewidth}
    \centering
    \includegraphics[width=0.75\linewidth]{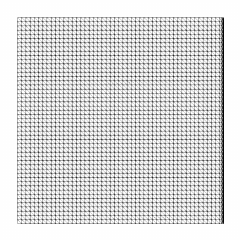}

    \caption{FE mesh, refined to capture shocks.}
  \end{subfigure}\hfill 
  \begin{subfigure}[t]{0.49\linewidth}
    \centering
    \includegraphics[width = \linewidth]{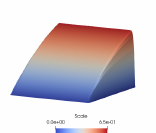}

    \caption{Steady state FOM solution at $t=20$ 
    for $\mathrm{Pe} =$ 6,000.}
  \end{subfigure}
    
    \caption{CDR test case: FOM mesh and steady state solution. FOM solution is rotated for visualization. 
   }
    \label{fig:CDR_FOM_Steady}
\end{figure}

\subsubsection{%
CDR Standard OpInf Setup}

To provide a more comprehensive assessment of EFR-OpInf, we consider a setting that combines parametric generalization with future-state prediction. 
Specifically, the system is parameterized by the diffusion coefficient $\epsilon$, and the ROM is evaluated both at an unseen parameter instance and for future-state predictions at that parameter value.
We train the ROM using $N_\epsilon = 4$ FOM trajectories corresponding to $\mathrm{Pe}\in \{5000,5500,6500,7000\}$ and evaluate its predictive performance at the unseen parameter value $\mathrm{Pe} = $ 6,000. 

We employ OpInf to learn a linear ROM of the form,
\begin{equation}
\frac{\mathrm{d} \state }{\mathrm{d} t} = \constantOP +  \epsilon \linearOP \state\, + \widehat{B} \state\,,
\label{eq:Opinf_dynamical_system_linear_CDR}
\end{equation}
where $\epsilon$ denotes the multiplicative inverse of the \peclet number. 
Our OpInf model thus consists of a parameter-dependent $\linearOP$ and a parameter-independent $\widehat{B}$, thereby isolating the effects of \peclet variation in $\linearOP$, while $\widehat{B}$ captures the contributions of terms involving the fixed parameters $\mathbf{b}$ and $\sigma$.
The operators defining~\eqref{eq:Opinf_dynamical_system_linear_CDR} 
are the solution to the minimization problem, 
\begin{equation}\label{eq:Opinf_CDR_Learning_Problem}
    \min_{\constantOP, \linearOP, \widehat{B}} \sum_{j = 1}^{4}\sum_{i=0}^{N_t} \left \|\frac{\mathrm{d} \oneDstate(t_i)  }{\mathrm{d} t}  - \constantOP - \epsilon_j \linearOP \oneDstate(t_i) -\widehat{B} \oneDstate(t_i)   \right \|_2^2 + \beta_1 (\| \linearOP  \|_F^2 + \| \widehat{B}\|_F^2 +  \| \constantOP \|_F^2)\,,
\end{equation}
where $\{\epsilon_j\}_{j=1}^4$ denote the $N_\epsilon=4$ training parameter instances. 
Similar formulations, extended to the multi-parameter case, are studied 
in~\cite{mcquarrie2023parametric} for a parametric heat equation and neuron model, as well as in \cite{vijaywargiya2025tensorparametrichamiltonianoperator} for a parametric heat equation and wave equation. 
The OpInf ROMs employed here are constructed using mean-centered snapshots because this transformation was found to improve model stability, although it is not required to enforce the boundary conditions. 
More generally, snapshot transformations are commonly used in OpInf to improve predictive performance or impose a desired model structure~\cite{Qian2019TranformandLearn,QIAN20201LiftandLearn}, and have been  
previously shown to enhance the stability of OpInf models across various complex applications such as solar-wind-stream prediction~\cite{ISSAN2023shiftedopinf} and rocket combustion simulations~\cite{farcas2023parametric,FARCAS2025109619,McQuarrie2021Regularized,QFW21}.

For comparison purposes, we will consider two standard OpInf scenarios: the first, unregularized OpInf, will set $\beta_1 = 0$. 
The second, strongly regularized OpInf, will find $\beta_1$ such that the OpInf error is minimized in the training region $(0,2]$ and the OpInf ROM remains stable until the end of the validation horizon, $t = 20.0$. 
Given these parameters, the strongly regularized $\beta_1$ was found to be $4.835$ which gave the smallest training error among the models %
satisfying the coefficient growth criterion
\eqref{eq:opinf_growth_ratio} over the validation interval, with
$\tau=1.2$ after searching 100 log-spaced values in the parameter range $[10^{-10}, 10^{3}]$ . 
All OpInf models will use $r = 10$, which captures $99.99\%$ of the energy in the POD modes, and are integrated in time using Runge--Kutta 4 scheme with a timestep of $\Delta t = 0.01$, identical to the FOM timestep. 

\subsubsection{CDR EFR-OpInf Setup}
\begin{figure}[!h]
    \centering
    \begin{subfigure}[t]{0.49\linewidth}
    \centering
    \includegraphics[width=0.75\linewidth]{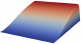}

    \caption{First POD mode}
  \end{subfigure}\hfill 
  \begin{subfigure}[t]{0.49\linewidth}
    \centering
    \includegraphics[width = 0.75\linewidth]{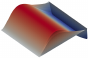}

    \caption{Second POD mode}
  \end{subfigure}
  \begin{subfigure}[t]{0.49\linewidth}
    \centering
    \includegraphics[width=0.75\linewidth]{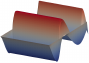}

    \caption{Fifth POD mode}
  \end{subfigure}\hfill 
  \begin{subfigure}[t]{0.49\linewidth}
    \centering
    \includegraphics[width = 0.75\linewidth]{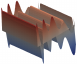}

    \caption{10th POD mode}
  \end{subfigure}
    
    \caption{CDR test case: POD modes. Color, scale, and orientation is individual to each basis function for clearer visualization.}
    \label{fig:CDR_POD_modes}
\end{figure}
The unregularized OpInf ROM serves as the baseline model for EFR-OpInf.
Additionally, EFR-OpInf-DF requires no parameter tuning, since the EFR parameters $\delta$ and $\chi$ are prescribed directly from theoretically motivated scalings.
Analysis of EFR formulations at the FOM level suggests the relaxation scaling $\chi=\mathcal{O}(\Delta t)$~\cite{neda2012_EFR_analysis}, while a standard choice for the filter radius is $\delta=\mathcal{O}(h)$~\cite{berselli2006LES,Sagaut2006LES}.
Accordingly, we set $\chi=\Delta t=0.01$ for all EFR-OpInf models and $\delta=h=0.02$ for the EFR-OpInf-DF and EFR-OpInf-HDF models, without performing any additional parameter search.
No filtering is applied over the training horizon, $t\leq2.0$, where the underlying OpInf ROM accurately reproduces the observed dynamics.
The prescribed EFR stabilization is activated only once the ROM enters the validation regime, where extrapolation errors and instabilities may develop rapidly.
Consequently, EFR-OpInf begins its future-state prediction from an accurate reduced state and introduces stabilization only where it is needed to support robust temporal extrapolation.

This leaves only one hyperparameter for EFR-OpInf: the reduced projection dimension $r'$ for the EFR-OpInf-Proj and EFR-OpInf-HDF models. 
Here we consider nine possible values $(r' < r = 10)$ whose solution is given by inspection of the POD basis plotted in Figure \ref{fig:CDR_POD_modes}. 
Inspection of the POD basis shows that all modes beyond the first exhibit pronounced oscillatory structure, giving the basis a Fourier-like character. 
The first mode is therefore the only one dominated by a smooth, non-oscillatory structure, which motivates our choice of $r'=1$ for both EFR-OpInf-Proj and EFR-OpInf-HDF. 

\subsubsection{%
CDR Comparison between OpInf and EFR-OpInf Models}
Figure~\ref{fig:CDR_errors} plots the $L^2$ errors over time for all considered OpInf models. 
Note that the errors for the unregularized OpInf model and all EFR-OpInf models are identical until the end of training data and diverge after the filter methods are turned on. 
The strongly regularized OpInf model is up to two orders of magnitude worse than all other OpInf models within the training region. 
In contrast, all EFR-OpInf models yield higher accuracy within the training region since they build on unregularized OpInf.  
Although strong regularization slows the error growth relative to the unregularized OpInf model, both predictions deteriorate and exhibit a persistent increase in error over time. 
On the other hand, all  EFR-OpInf models are stable and reach a steady state in terms of error by around $t = 5.0$. 
Importantly, their errors do not grow in time.  

The EFR-OpInf-DF model performs similarly well to EFR-OpInf-Proj, with
EFR-OpInf-HDF performing the worst of the EFR-OpInf models. %
However, all EFR-OpInf models are stable and display no more than $10^{-4}$ relative $L^2$ error throughout the duration of the simulation, far better than either standard OpInf method is able to perform.
While EFR-OpInf-HDF remains stable and is more accurate than the strongly regularized OpInf model, it is outperformed by the other filtered ROMs for this CDR problem. 
This result indicates that the relative effectiveness of the different filtering strategies can depend on the underlying problem and flow features. 
As shown in Section~\ref{sec:NSE_FPC_Compare}, EFR-OpInf-HDF exhibits a better performance for the NSE equations, where it becomes competitive with the best-performing filtered ROMs.
\begin{figure}[htt]
    \centering
        \includegraphics[width=\linewidth]{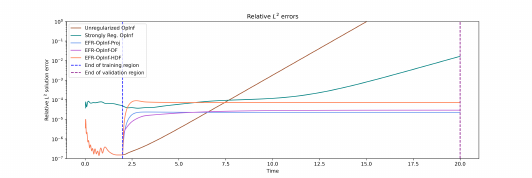}
        \caption{CDR test case: comparison of OpInf models at P\'eclet number 6,000 ($r = 10$).}
    \label{fig:CDR_errors}
\end{figure}

Figure \ref{fig:CDR_Solution_Compare} contrasts the solutions of the EFR-OpInf-DF model and the best performing standard OpInf model, i.e., the strongly regularized OpInf model.   
While the EFR-OpInf-DF model remains stable and displays no oscillations, the standard OpInf model displays oscillations as early as $t = 5$, including sharp oscillations near the boundary layer by $t = 10$ and complete breakdown by $t = 20$. 
Recall in addition that the stable and accurate EFR-OpInf-DF requires no tuning: there is no regularization baked into the underlying OpInf problem, and each of the EFR parameters $\delta$ and $\chi$ are chosen according to scalings obtained from theoretical error bounds without any parameter search. 
\begin{figure}[ht]
  \centering

  \begin{subfigure}{\linewidth}
    \centering
    \includegraphics[width=0.75\linewidth]{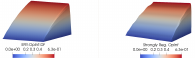}
    \caption{EFR-OpInf-DF (left) and Strongly Regularized OpInf solutions (right) at $t = 10$}
  \end{subfigure}\par\medskip

  \begin{subfigure}{\linewidth}
    \centering
    \includegraphics[width=0.75\linewidth]{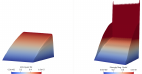}
    \caption{EFR-OpInf-DF (left) and Strongly Regularized OpInf solutions (right) at $t = 20$}
  \end{subfigure}
  \caption{CDR test case: OpInf solutions at P\'eclet number 6,000. Compare with Figure \ref{fig:CDR_FOM_Steady}.}
    \label{fig:CDR_Solution_Compare}
  \end{figure}

\subsection{Navier--Stokes Equations: 2D Flow Past a Cylinder %
in the Transitional Regime} %
    \label{sec:NSE}

We next test the proposed EFR-OpInf method on the more complex setting of the unsteady incompressible Navier--Stokes equations. 
Specifically, we assess its ability to construct predictive ROMs for flow past a cylinder, a canonical benchmark problem widely used in the ROM literature~\cite{gunzburger2019EFR_ROM,Mohebujjaman2019physically,strazzullo2021_EFRConsistency,Wells2017EF}.

\subsubsection{FPC FOM Setup}

The training data are derived from a FE simulation using Taylor--Hood elements.
Figure~\ref{fig:cylinder_mesh} plots the underlying mesh. 
The discrete velocity field has 8,024 DoFs per velocity component, yielding a total of 16,048 DoFs for the full-order solution.
The problem setup closely follows a FEniCSx implementation tutorial that 
can be found online at~\cite{DokkenFPC}.
The length of the domain is 2.2, the height is 0.41, and the cylinder has a radius of 0.05 with a center point of (0.2,0.2), considering the bottom-left of the rectangular domain to be the origin.

\begin{figure}[httbp]
  \centering

  \begin{subfigure}{0.48\linewidth}
    \centering
    \includegraphics[width=\linewidth]{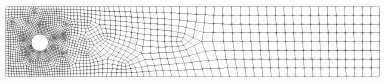}
    \caption{FE Mesh}
    \label{fig:cylinder_mesh}
  \end{subfigure}\hfill
  \begin{subfigure}{0.48\linewidth}
    \centering
    \includegraphics[width=\linewidth]{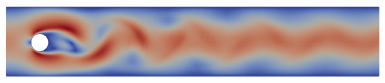}
    \caption{FOM velocity magnitude at time $t = 8.0$}
    \label{fig:Cylinder_FOM_sol}
  \end{subfigure}

  \caption{NSE FPC: Finite element mesh (left) and full-order model velocity magnitude at $t = 8.0$ (right) with magnitude scale [0,2.2].}
  \label{fig:NSE_mesh_and_FOM}
\end{figure}

The simulation is run from a cold start with a parabolic inflow condition for $u_x,$ the $x$-component of $\bu$:
$$
u_x(t, y) =
\begin{cases}  \dfrac{6 \sin(\pi t /8) y(y-0.41)}{0.41^2}, \;\; t < 4, \\ 1.5, \;\; t \geq 4,
\end{cases}
$$
where the $y$-component, $u_y,$ is set equal to zero. 
This provides a maximum flow velocity of $1.5$ in the horizontal direction, resulting in a Reynolds number of $\mathrm{Re} = 100$. 
The FOM is solved with a pressure-correcting method using a Crank--Nicolson scheme with a FOM $\Delta t$ of $0.0005$ over $t \in (0.0, 20.0]$.
The mesh is generated according to a scheme where the mesh elements near the cylinder have a target mesh resolution of $h = 0.05/3 \approx 0.02$, one-third of the size of the cylinder radius. 
The snapshots are uniformly spaced in time with $\Delta t=0.005$, which is also used as the timestep for all OpInf and EFR-OpInf ROM simulations. The FOM takes approximately 15 minutes to run all 20 seconds using 4 MPI processes on a laptop.

We collect snapshots over the time span of $[4.0, 20.0]$, partitioned into the training interval $[4.0, 6.0)$, the validation interval $[6.0, 10.0)$, and
the additional testing interval $[10.0, 20.0]$. The FOM produces 40,000 snapshots over the interval $[0.0, 20.0]$, but we retain every tenth snapshot to reduce storage requirements. This results in 400 training snapshots over $[4.0, 6.0)$, 800 validation snapshots over $[6.0, 10.0)$, and 2,000 testing snapshots over $[10.0, 20.0]$.

\subsubsection{FPC Standard
OpInf Setup}
\label{sec:NSE_OpInf_Setup}

We retain the first $r=20$ POD modes, which capture 99.6\% of the total snapshot energy.
The magnitudes of the $1^{\rm st}$, $5^{\rm th}$, $10^{\rm th}$, and $15^{\rm th}$ POD velocity basis vectors are plotted in Figure~\ref{fig:rom-basis-funcs}. %
As expected, the low-index POD modes are dominated by large-scale flow structures, whereas the higher-index modes display progressively more oscillatory, small-scale features. The NSE OpInf model is constructed according to the process discussed in Section~\ref{subsec:OpInf_summary}. 

\begin{figure}[httbp]
  \centering
  \begin{subfigure}{0.48\linewidth}
    \centering
    \includegraphics[width=\linewidth]{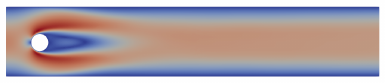}
    \caption{Mean (non-fluctuating component)}
  \end{subfigure}

  \vspace{0.8em} %

  \begin{subfigure}{0.48\linewidth}
    \centering
    \includegraphics[width=\linewidth]{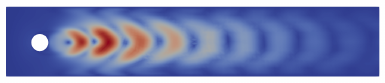}
    \caption{1st ROM POD Mode}
  \end{subfigure}\hfill
  \begin{subfigure}{0.48\linewidth}
    \centering
    \includegraphics[width=\linewidth]{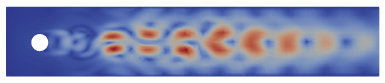}
    \caption{5th ROM POD Mode}
  \end{subfigure}

  \vspace{0.8em} %

  \begin{subfigure}{0.48\linewidth}
    \centering
    \includegraphics[width=\linewidth]{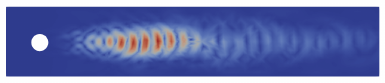}
    \caption{10th ROM POD Mode}
  \end{subfigure}\hfill
  \begin{subfigure}{0.48\linewidth}
    \centering
    \includegraphics[width=\linewidth]{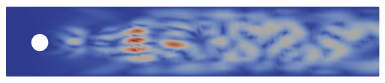}
    \caption{15th ROM POD Mode}
  \end{subfigure}

  \caption{NSE FPC: Magnitudes of ROM POD modes and mean of snapshot data. Basis functions are derived from the fluctuating component of flow and color map scaling varies by function.}
  \label{fig:rom-basis-funcs}
\end{figure}

We will consider two OpInf regularization strategies for $\hyperparameter$: the weakly regularized OpInf method (model~\ref{model:weakly_regularized}) and the strongly regularized OpInf method (model~\ref{model:strongly_regularized}). The regularization parameters are obtained from a grid search: recalling that $\hyperparameter = [\beta_1, \beta_2]^T$ where $\beta_1$ regularizes $\constantOP$ and $\linearOP$, and $\beta_2$ regularizes $\quadOP$, we take 10 log-spaced values between $(10^{-5}, 10^2)$ for $\beta_1$ and 15 log-spaced values between $(10^{-2}, 10^4)$ for $\beta_2$, choosing the optimal $\beta_1, \beta_2$ values in this range. After this grid search, $\hyperparameter = [0.1, 4.6]^T$ was chosen for the weakly regularized EFR-OpInf, and $\hyperparameter = [31, 464]^T$ was found which gives the smallest training error among the models
satisfying the coefficient growth criterion
\eqref{eq:opinf_growth_ratio} over the validation interval, with
$\tau=1.2$ for strongly regularized OpInf.

The solutions of these models at time $t = 8.0$, midway through the validation region, are plotted in Figure \ref{fig:NSE_cylinder_opinf_sols}. Comparing this against the FOM solution at the same timestep in Figure \ref{fig:Cylinder_FOM_sol}, it is clear that neither OpInf solution is physical. The differing effects of regularization %
are also evident: the strongly regularized OpInf ROM displays considerably lower velocity magnitudes as compared to the unregularized OpInf ROM. However, in this case of limited training data, neither OpInf ROM properly learns the dynamics of the flow past a cylinder problem, even within the bounds of the validation region. KE and $L^2$ error results will follow in Figures \ref{fig:FPC_EFR_OpInf_KE} and \ref{fig:FPC_EFR_OPINF_L2}, but each OpInf model displays unbounded growth in both metrics over the full test region, with the weakly regularized OpInf model diverging within the validation region. 

\begin{remark}[ROMs for Transitional Flows]
The physical and numerical parameters chosen (such as the inflow condition from cold start) result in a FOM velocity field which is still developing in the training regime of $t \in [4.0, 6.0]$ and is not statistically steady (see FOM KE in Figure \ref{fig:FPC_EFR_OpInf_KE}). 
This is not the traditional setting for ROMs of FPC settings (\cite{gunzburger2019EFR_ROM,Mohebujjaman2019physically,strazzullo2021_EFRConsistency,Wells2017EF}), where snapshots are collected within the statistically steady regime and used to predict data from the same statistical regime. 
This is fundamentally different from the setting we consider here, where snapshots collected from the transitional regime are used to predict into the statistically steady region. It is promising that EFR-OpInf can provide significant performance improvements for this more challenging than usual FPC problem. 
\end{remark}

\begin{figure}[httbp]
  \centering
  \begin{subfigure}[c]{0.15\textwidth}
    \centering
    \includegraphics[width=\textwidth, height=0.25\textheight, keepaspectratio]{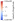}

  \end{subfigure}%
  \hfill
  \begin{subfigure}[c]{0.82\textwidth}
    \centering
    \begin{subfigure}[b]{\textwidth}
      \centering
      \includegraphics[width=0.8\textwidth]{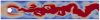}
      \caption{Weakly Regularized OpInf}
      \label{fig:opinf_weak_sol}
    \end{subfigure}
    
    \vspace{0.5em} %
    
    \begin{subfigure}[b]{\textwidth}
      \centering
      \includegraphics[width=0.8\textwidth]{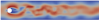}
      \caption{Strongly Regularized OpInf}
      \label{fig:opinf_strong_sol}
    \end{subfigure}
  \end{subfigure}

  \caption{NSE FPC: OpInf solutions at $t = 8.0$. Compare to Figure \ref{fig:Cylinder_FOM_sol}.}
  \label{fig:NSE_cylinder_opinf_sols}
\end{figure}

\subsubsection{FPC EFR-OpInf Setup}

We now discuss the choice of parameters for each of our EFR models. For EFR-OpInf-Proj, an inspection of the basis (Figure \ref{fig:rom-basis-funcs}) suggests that $r' = 10$ (half of ROM dimension $r = 20$) is a reasonable transition point between smoother and more oscillatory basis components. For EFR-OpInf-DF, we leverage the numerical analysis available at the FOM level \cite{berselli2006LES} and choose $\delta = h \approx 0.02$. The main parameter that we will allow to vary is $\chi$, which controls the amount of filtering used in the relaxation step, with smaller $\chi$ corresponding to less filtering. To maintain a consistent comparison with standard OpInf regularization methods, we will use a similar ROM coefficient-based strategy as to that which was described for standard OpInf in \ref{sec:opinf_regularization}. 

To accomplish this, we select 15 log-spaced values of $\chi$ in the range of $[0.05, 1.0]$, and run the EFR-OpInf ROM over the validation region. We compute the growth of the online EFR-OpInf ROM and compare it to growth seen within training according to Equation~\eqref{eq:opinf_growth_ratio}. Selecting a maximum growth ratio of $M = 1.2$, we discard any EFR parameter choices that result in excessive growth in the validation region. From those parameter values remaining, we choose $\chi$ such that the corresponding EFR-OpInf ROM shows the least growth within the validation region. After running this process for each of our methods, we obtain the following values of $\chi$: for EFR-OpInf-Proj, $\chi = 0.1805$, for EFR-OpInf-DF, $\chi = 0.05$, and for EFR-OpInf-HDF, $\chi = 0.05$. These parameter choices are summarized in Table \ref{tab:FPC_EFR_Parameters}.
\begin{table}[ht]
    \centering
    \begin{tabular}{c|c|c|c}
        & $r'$ & $\delta$ & $\chi$ \\ \hline
        EFR-OpInf-Proj & 10 & - & 0.1805 \\
        EFR-OpInf-DF & - & 0.02 & 0.05 \\
        EFR-OpInf-HDF & 10 & 0.02 & 0.05
    \end{tabular}
    \caption{NSE FPC: EFR-OpInf Parameter Values.}
    \label{tab:FPC_EFR_Parameters}
\end{table}

\subsubsection{FPC Comparison between OpInf and EFR-OpInf Models}
\label{sec:NSE_FPC_Compare}

Here we compare the accuracy, stability, and online costs of our various OpInf and EFR-OpInf models for the flow past a cylinder problem. In Figure \ref{fig:FPC_EFR_OPINF_SOLS}, we compare each of our EFR-OpInf models against the FOM solution at $t = 15.0$, well into the testing region. Each of the models is stable, but clearly EFR-OpInf-DF is overly diffusive. In contrast, both EFR-OpInf-Proj and EFR-OpInf-HDF 
remain visually close to the FOM solution 9 seconds after leaving the training data. 

\begin{figure}[httbp]
  \centering
  \begin{subfigure}[c]{0.1\textwidth}
    \centering
    \includegraphics[width=\textwidth, height=0.45\textheight, keepaspectratio]{FPC_figure_8_scale.pdf}
  \end{subfigure}%
  \hfill
  \begin{subfigure}[c]{0.90\textwidth}
    \centering
    
    \begin{subfigure}[b]{0.48\textwidth}
      \centering
      \includegraphics[width=\textwidth]{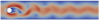}
      \caption{FOM Solution}
    \end{subfigure}\hfill
    \begin{subfigure}[b]{0.48\textwidth}
      \centering
      \includegraphics[width=\textwidth]{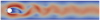}
      \caption{EFR-OpInf-Proj}
    \end{subfigure}

    \vspace{0.8em} %

    \begin{subfigure}[b]{0.48\textwidth}
      \centering
      \includegraphics[width=\textwidth]{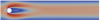}
      \caption{EFR-OpInf-DF}
    \end{subfigure}\hfill
    \begin{subfigure}[b]{0.48\textwidth}
      \centering
      \includegraphics[width=\textwidth]{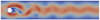}
      \caption{EFR-OpInf-HDF}
    \end{subfigure}
  \end{subfigure}
  \caption{NSE FPC: Comparison EFR-OpInf Velocity Magnitudes to FOM at $t = 15.0$.}
  \label{fig:FPC_EFR_OPINF_SOLS}
\end{figure}

Comparing KEs in Figure \ref{fig:FPC_EFR_OpInf_KE}, we see that  weakly regularized OpInf provides an excellent match to the FOM KE within the training window but diverges quickly past training, while strongly regularized OpInf provides lower accuracy within the training window and ultimately fails to remain stable at the testing horizon. By using weakly regularized OpInf as a springboard for our EFR models, each of our EFR methods is able to retain optimal KE within training window while remaining long-term stable. For these parameter values, EFR-OpInf-HDF and EFR-OpInf-Proj provide %
comparable solutions, with EFR-OpInf-HDF coming out on top, while EFR-OpInf-DF is overly diffusive. 

\begin{figure}
    \centering
    \includegraphics[width=\linewidth]{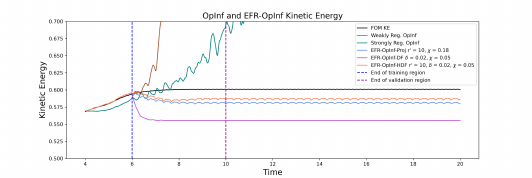}
    \caption{NSE FPC: OpInf, EFR-OpInf, and FOM kinetic energies.}
    \label{fig:FPC_EFR_OpInf_KE}
\end{figure}

The story for $L^2$ errors in Figure \ref{fig:FPC_EFR_OPINF_L2} is 
similar to that for KE, with weakly regularized OpInf spiking to 40\% error by the end of training window, as the excessively high $\hyperparameter$ damages the solution accuracy. In contrast, the error curves for all ERF-OpInf models remain flat over a greatly extended time horizon. In cases where it is difficult for standard OpInf regularization techniques to balance stability and accuracy, adding this extra layer of regularization lets OpInf focus on accuracy while handling the stability through filtering. 

\begin{figure}
    \centering
    \includegraphics[width=\linewidth]{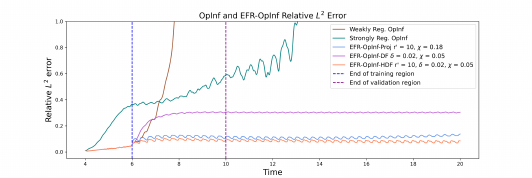}
    \caption{NSE FPC: OpInf and EFR-OpInf relative velocity $L^2$ errors.}
    \label{fig:FPC_EFR_OPINF_L2}
\end{figure}

Our last comparison, Figure \ref{fig:Pareto_Vary_Chi} is a Pareto plot that gives an idea of model performance and online time cost. 
Variations in performance are 
created by altering the relaxation parameter $\chi$ 
shared across all EFR-OpInf models, with
minimum $\chi = 0.005 = \Delta t$ and 
maximum value of $\chi = 1.0$ (no relaxation). 
In this figure, the `OpInf baseline' label refers to the strongly regularized, best performing standard OpInf model. We will compare the models in the validation region, as this gives the fairest comparison to standard OpInf techniques in the low-data setting. 

In Figure \ref{fig:Pareto_Vary_Chi}, both KE and $L^2$ relative errors are computed as averages over the time horizon of $[4.0, 10.0]$, the training and validation regions. This provides a comparison in the regions where OpInf can be expected to perform the best.
In Figure \ref{fig:Pareto_Vary_Chi}, the wall-clock time refers to the online cost of time-integrating the ROM from $t = 4.0$ to $t = 10.0$, which is 1,200 timesteps of RK4. For these choices, the computed quantities were finite (i.e., the model did not diverge) for all but the first $\chi$ value for EFR-OpInf-Proj, and all but the first two $\chi$ values for both differential filters. In general, we observed that the EFR-OpInf models were consistently stable after a certain threshold of $\chi$ was reached, which for this problem was larger than $\Delta t$. Previous work on filtering for ROMs has suggested~\cite{Moore2025_ADLROM} that ROMs can tolerate, and sometimes require, %
more filtering than corresponding FOM filtering techniques.

Observe that the time costs are extremely similar for each of the OpInf and EFR-OpInf models considered: while the errors vary 
more than an order of magnitude, the runtimes for every OpInf model, EFR or not, are within 0.08 seconds of each other. 
EFR-OpInf-Proj is the simplest filter and nearly identical to the baseline OpInf model at 0.20-0.21 seconds, but both EFR-OpInf-DF and EFR-OpInf-HDF average out to 0.23-0.24 seconds online.
This shows that there is very little online cost separating the projection and differential filters. %

For the kinetic energy errors in Figure \ref{fig:Pareto_Vary_Chi}(a), only the overly dissipative EFR-OpInf-DF model displays worse 
kinetic energy error 
than the OpInf model. 
EFR-OpInf-Proj provides reasonable KE improvements over standard OpInf, and EFR-OpInf-HDF provides the best KE matching of all the models considered. 
Neither the EFR-OpInf-Proj nor EFR-OpInf-DF models display significant sensitivity to $\chi$. 
The performance of the EFR-OpInf-HDF model is similar to the 
EFR-OpInf-Proj model for low $\chi$, and it can be tuned to perform much better (in terms of KE) than all competing models for 
values of $\chi$ approaching 1. Since $\chi=1$ corresponds to using the filtered state without relaxation, these results suggest that the relaxation parameter may be eliminated from the EFR-OpInf-HDF model if the goal is KE-optimal models. 
This would remove the principal hyperparameter tuning required for EFR-OpInf-HDF, representing an advantage over both the projection and differential filters. 

For the $L^2$ errors in Figure \ref{fig:Pareto_Vary_Chi}(b), each of the EFR-OpInf models 
outperform the standard OpInf model, with the EFR-OpInf-Proj and EFR-OpInf-HDF 
outperforming standard OpInf by an order of magnitude or more for most values of $\chi$. While the EFR-OpInf-Proj model displays even less sensitivity to $\chi$ in terms of $L^2$ errors, the EFR-OpInf-HDF model displays a similar sensitivity to that 
seen in the kinetic energy errors, with all $\chi$ values tested superior to standard OpInf but some choices of $\chi$ better than others. In particular, the reverse of the kinetic energy case holds: smaller values of $\chi$ tend to produce better $L^2$ error values than larger values. While we exclude additional visualizations for space, 
applying EFR-OpInf-HDF with large $\chi$ 
produces a flow that is also representative of the FPC problem.  However, the EFR-OpInf-HDF velocity field is phase-shifted from the FOM solution, which causes pointwise $L^2$ errors to degrade.

\begin{figure}[htt]
  \centering

    \begin{subfigure}[t]{0.49\linewidth}
    \centering
    \includegraphics[width=\linewidth]{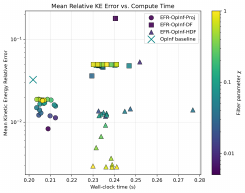}

    \caption{KE error vs online ROM time}
  \end{subfigure}\hfill 
  \begin{subfigure}[t]{0.49\linewidth}
    \centering
    \includegraphics[width=\linewidth]{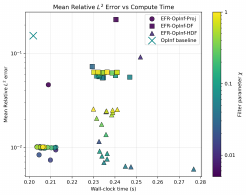}
    \caption{$L^2$ error vs online ROM time}
  \end{subfigure}
  \caption{NSE FPC: Pareto plots for mean relative $L^2$ and KE errors vs. online ROM cost for $t \in [4.0, 10.0]$ (training + validation region). OpInf baseline corresponds to the strongly regularized OpInf model. Models are colored by varying $\chi \in [0.005, 1.0]$. 
  }
  \label{fig:Pareto_Vary_Chi}
\end{figure}

The results of the runs for the Pareto chart in Figure \ref{fig:Pareto_Vary_Chi} are summarized in Table \ref{tab:pareto_table}, reporting the averages of both $L^2$ and KE errors as well as the minimums, and the $\chi$ values for which they occurred. 

\begin{table}[ h]
    \centering
    \begin{tabular}{l|c|c|c|c}
      Model   & Avg. $L^2$ Error & Avg. KE Error & Min. $L^2$ Error & Min. $KE$ Error \\
      \hline
      Strongly Reg. OpInf   & 15.43\% & 3.3\% & 15.43\% & 3.3\% \\
      \hline 
      EFR-OpInf-Proj & 1.15\% & 1.5\% & 0.73\% ($\chi = 0.015$) & 0.83\% ($\chi = 0.005$) \\ \hline 
      EFR-OpInf-DF & 6.07\% & 4.81\% & 5.48\% ($\chi = 0.05$) & 3.62\% ($\chi = 0.03$)  \\ \hline 
        EFR-OpInf-HDF & 1.92\% & 1.11\%  & 0.58\% ($\chi = 0.04$) & 0.28\% ($\chi = 0.8$)
    \end{tabular}
    \caption{NSE FPC: Summary of time-averaged errors over validation region from Figure \ref{fig:Pareto_Vary_Chi}. }
    \label{tab:pareto_table}
\end{table}

\subsection{Navier--Stokes Equation: 3D Minimal Channel Flow}
\label{sec:NSE_mfu}

We next consider the 3D minimal channel flow scenario at $\rm Re=5000$, which exhibits strongly turbulent dynamics while retaining a computationally tractable domain size. 
Minimal channel flow has been widely used as a reduced computational setting for studying near-wall turbulence, since it preserves important turbulent structures while significantly lowering the computational cost compared with a full channel flow simulation~\cite{jimenez1991minimal}.

\subsubsection{MCF FOM Setup}
Following the setup in~\cite{jimenez1991minimal}, the streamwise and spanwise lengths of the channel are set to $0.6\pi$ and $0.18\pi$, respectively, and the channel half-height is set to $1$. 
The boundary conditions are periodic in the streamwise and spanwise directions, with no-slip boundary conditions imposed at the channel walls. %
For this scenario, time is nondimensionalized by the convective time
scale $t_c=L/U$, where $L$ and $U$ are the characteristic length and velocity
scales, respectively. 
Thus, the nondimensional time $t^\star=t/t_c=tU/L$ is measured in convective time units (CTUs).

The FOM simulation is carried out using nekRS~\cite{fischer2022nekrs}. 
The spatial discretization uses $576$ spectral elements with polynomial order
$N=9$, resulting in approximately $5.76\times 10^5$ spatial degrees of
freedom per velocity component. 
Thus, each three-component velocity snapshot contains approximately $1.73\times 10^6$ degrees of freedom.
The FOM uses a semi-implicit BDF$k$/EXT$k$ time discretization with timestep $\Delta t = 0.0025$ and orders up to $k=3$, where BDF$k$ is used for the time derivative, EXT$k$ for the advection and forcing terms, and an implicit treatment for the
dissipation terms.

We focus on the time interval $[3000, 3100]$ CTUs, after the solution has reached a statistically steady regime.  
The ROM is set up so that $\bu^r(t,\mathbf{x}) = \umean + \sum_{i=1}^r a_i(t)\basisfunc_i(\mathbf{x}),$
where the basis functions $\basisfunc_i$ are obtained through POD on the snapshot matrix of the fluctuating component of the high-dimensional FOM data. 
We divide the full time interval into the
training interval $[3000,3025]$, the validation interval $(3025,3050]$, and the
testing interval $(3050,3100]$. Over the training interval, snapshots are sampled every $0.0125$ CTUs, leading to a total of $2001$ snapshots, including the initial condition at $t=3000$. 
These snapshots are used to construct the POD basis and train both the weakly and strongly regularized OpInf models. 
The FOM simulation is then continued through the validation interval, from which one representative snapshot is retained to guide the selection of the EFR-OpInf filtering parameter $r'$. 
Finally, the testing interval is used to evaluate the long-time predictive accuracy and stability of the resulting OpInf and EFR-OpInf models.

\subsubsection{MCF Standard OpInf Setup}
\label{sec:MFU_standard_OpInf}

The POD basis is constructed with a ROM dimension of $r=29$, which retains approximately $60\%$ of the cumulative energy of the available modes. 
In our preliminary numerical investigation, we have also tested the ROM dimensions $r=44$ (which retains approximately $70\%$ of the cumulative energy) and $r=67$ (which retains about $80\%$ of the cumulative energy).
Since the resulting numerical behavior was qualitatively consistent across these choices, we present only the results for $r=29$ for clarity.
This choice places the ROM in a more strongly under-resolved regime, representative of the challenges encountered when modeling turbulent flows with low-dimensional approximations.
We note that the retained energy for $r=29$ is substantially lower than in the %
FPC example, reflecting the richer dynamics and broader range of active scales present in the minimal channel flow.
Consequently, this example provides a particularly challenging setting for assessing the stability and predictive performance of the OpInf and EFR-OpInf models.

For the standard OpInf models, we consider both the weakly and strongly regularized versions described in Section~\ref{sec:opinf_regularization}.  
The regularization parameters are obtained through a grid search.  
Recalling that $\hyperparameter = [\beta_1, \beta_2]^T$, where $\beta_1$ regularizes $\constantOP$ and $\linearOP$, and $\beta_2$ regularizes $\quadOP$, we take 10 logarithmically spaced values between $10^{-5}$ and $10^2$ for $\beta_1$, and 15 logarithmically spaced values between $10^{-5}$ and $10^2$ for $\beta_2$.
After this grid search, the weakly regularized model uses
$\hyperparameter = [3.59 \times 10^{-4}, 10^{-3}]^T$, which gives the smallest training error over the grid. 
The strongly regularized model uses $\hyperparameter = [10^{-5}, 10^{-2}]^T$, which gives the smallest training error among the models %
satisfying the coefficient growth criterion
\eqref{eq:opinf_growth_ratio} over the validation interval, with
$\tau=1.2$.
We note that the selected value of
$\beta_1$ lies at the lower bound of the search interval. 
Therefore, this parameter choice is optimal only over the grid considered, and extending the search interval could yield a different strongly regularized model. 

Figure~\ref{fig:Reg_OpInf_mfu_ke_compare_proj} compares the kinetic energy of the FOM, the weakly regularized, and the strongly regularized OpInf models. 
In the training region, the weakly and strongly regularized OpInf models exhibit similar kinetic energy behavior, but both show a gap from the FOM energy due to the approximation error associated with $r=29$.
After the training window, the weakly regularized OpInf model overfits and quickly becomes unstable, whereas the strongly regularized OpInf model remains stable throughout the validation region but eventually becomes unstable in the testing region. 

\begin{figure}[htt]
    \centering
    \includegraphics[width = \linewidth]{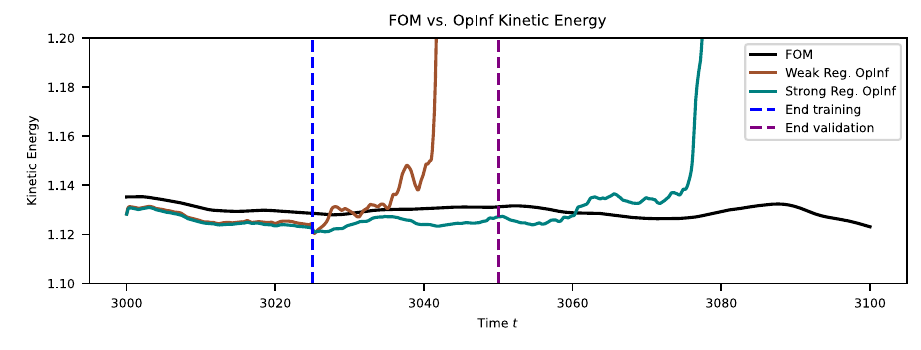}
    \caption{3D MCF at $\rm Re=5000$: Kinetic energy behavior of the FOM, the weakly regularized OpInf model ($\beta_1 = 3.59\times 10^{-4}, \beta_2 = 10^{-3}$), and the strongly regularized OpInf model ($\beta_1 = 10^{-5}, \beta_2 = 10^{-2}$).}
    \label{fig:Reg_OpInf_mfu_ke_compare_proj}
\end{figure}

Figure~\ref{fig:opinf_velocity_fields_t_3042} shows the velocity magnitude fields at $t=3042$ in the validation region. 
At this time instant, the weakly regularized OpInf prediction has already become non-physical, whereas the strongly regularized OpInf model still produces a reasonable velocity field.
\begin{figure}[htt]
    \centering
    \includegraphics[width=1.0\linewidth]{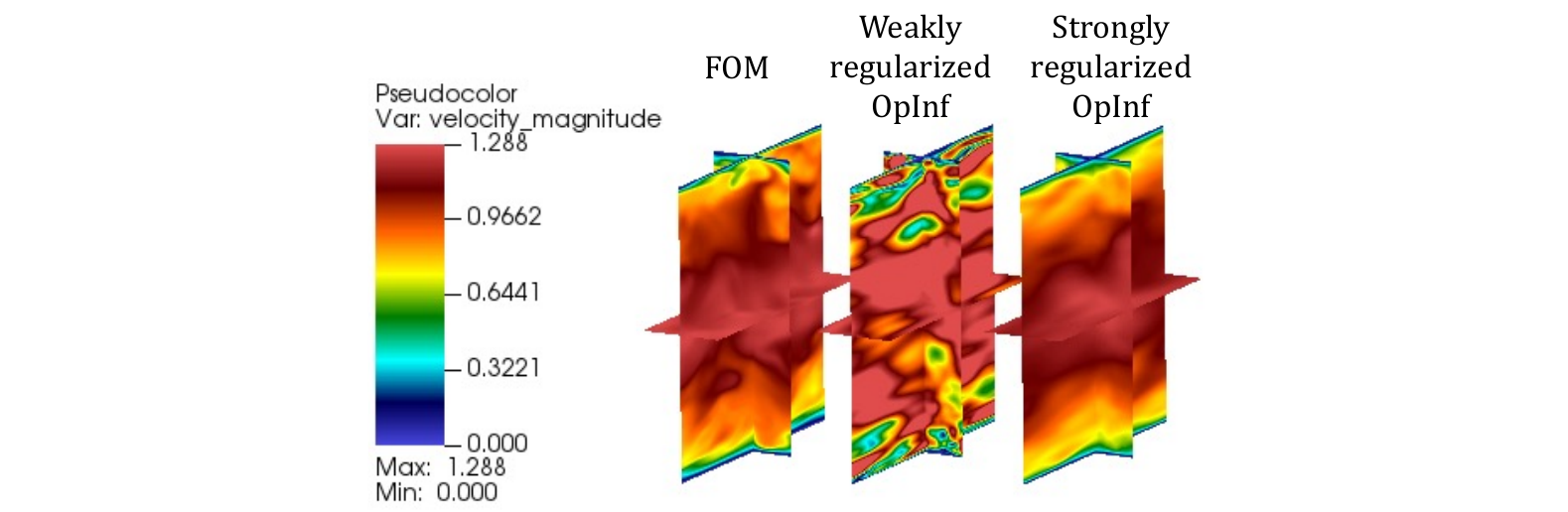}
    \caption{3D MCF at $\rm Re=5000$: Velocity magnitude fields computed by the FOM and predicted by the weakly and strongly regularized OpInf models at $t=3042$.}
  \label{fig:opinf_velocity_fields_t_3042}
\end{figure}

\subsubsection{MCF EFR-OpInf Setup}
\label{sec:MFU_EFR_OpInf_Setup}

For all three EFR-OpInf experiments, the underlying OpInf model is the weakly regularized model with
$\hyperparameter = [3.59 \times 10^{-4},\, 10^{-3}]^T$.
Online EFR stabilization is activated only after the ROM exits the training interval, $[3000,3025)$. 
Thus, the unfiltered OpInf ROM is evolved over the training interval, whereas EFR stabilization is applied only afterwards.
For EFR-OpInf-Proj, we use the ROM projection filter described in Section~\ref{sec:Proj_Filter}. 
For this scenario, we choose $r' = 11$. 
The cutoff dimension $r'$ is selected using the same strategy as in the previous scenario. 
We first consider a representative FOM velocity field at $t=3037.5$ in the validation region, project it onto the $r$-dimensional POD space, and examine how different cutoff values $r'$ modify the projected velocity field. 
The values of $r'$ for which visible changes occur are retained as candidates, giving $r' \in \{29,28,19,18,16,12,11,4,3\}$. 
For each candidate cutoff dimension, we vary the relaxation parameter $\chi$ over values proportional to the FOM timestep $\Delta t$. For each EFR-OpInf-Proj parameter pair
$(r',\chi)$, we compute the coefficient growth ratio defined in
\eqref{eq:opinf_growth_ratio} over the validation interval and retain the pair as admissible if this ratio does not exceed the prescribed tolerance $\tau=1.2$.
This selection process does not use kinetic energy information and is consistent with the selection strategy used for the strongly regularized OpInf baseline.
The admissible cutoff values are $r'=11$, $r'=4$, and $r'=3$. 
We choose the largest admissible cutoff, $r'=11$, together with $\chi=0.000625=0.25\Delta t$, so that the filter retains as many resolved POD modes as possible.

For EFR-OpInf-DF, we use the ROM differential filter defined in \eqref{eq:ROM_DF_Filter}. 
The filter radius is chosen as $\delta = 0.077$, which is half of the average mesh size of the spectral elements used in the FOM simulation. 
To keep the EFR-OpInf-DF setup consistent with EFR-OpInf-Proj, we use the same relaxation parameter, $\chi=0.000625=0.25\Delta t$. 
For EFR-OpInf-HDF, we use the hybrid ROM projection-differential filter introduced in Section~\ref{sec:hybrid_proj_diff_filter}. 
This filter applies the differential filter only to the high-index portion of the reduced state, namely the last $r-r'$ modes, while leaving the first $r'$ modes unchanged. 
For consistency with the projection and differential filter choices above, we set $r'=11$, $\delta=0.077$, and $\chi=0.000625=0.25\Delta t$. 

\subsubsection{MCF Comparison between OpInf and EFR-OpInf Models}
\label{sec:mfu_compare}

In this section, we compare the predictive performance and online costs of the standard OpInf models and the EFR-OpInf models for the minimal channel flow scenario. 
The comparison includes the weakly regularized OpInf model, the strongly regularized OpInf model, and the three EFR-OpInf models constructed using the ROM projection, ROM differential, and hybrid ROM projection-differential filters. 
We evaluate these models through the end of the testing interval and compare them in terms of kinetic
energy, relative velocity $L^2$ error, representative velocity fields, and
online computational cost.

Figure~\ref{fig:EFR_OpInf_KE_compare_mfu} compares the kinetic energy of the FOM, the OpInf models, and the three EFR-OpInf models. 
\begin{figure}[!h]
    \centering
    \includegraphics[width=\linewidth]{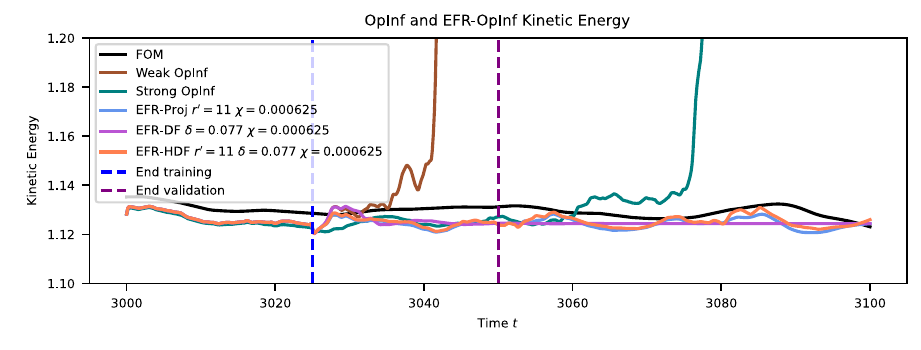}
    \caption{3D MCF at $\rm Re=5000$: Kinetic energy behavior of the FOM, the OpInf, and the EFR-OpInf models.}
    \label{fig:EFR_OpInf_KE_compare_mfu}
\end{figure}
As noted above, the weakly regularized OpInf model becomes unstable after the training region, while the strongly regularized OpInf model remains stable through validation but eventually diverges in the testing region. In contrast, all three EFR-OpInf models remain stable over the full interval $[3000,3100]$. This demonstrates that applying EFR stabilization during the online stage can improve the long-time stability of the weakly regularized OpInf model without modifying the learned operators. 
Among the EFR-OpInf models, EFR-OpInf-Proj and EFR-OpInf-HDF perform similarly, with kinetic energies tracking the FOM much closer than EFR-OpInf-DF. 
These numerical results suggest that, for this $\chi$ choice, the differential filter in the EFR-OpInf-HDF has a negligible effect. 
Although the kinetic energy of EFR-OpInf-DF remains bounded, the model is overly dissipative, and its energy eventually settles to a constant value. The %
hybrid ROM projection-differential filter mitigates this over-dissipation by applying differential filtering only to the high-index modes.

The relative $L^2$ errors for velocity, plotted in
Figure~\ref{fig:EFR_OpInf_L2error_compare_mfu}, show a consistent trend. 
In the training region, all models perform similarly, with relative errors around $6\%$--$8\%$.
The errors rise to roughly $8\%$--$12\%$ upon entering the validation region, likely because the ROM subspace constructed from training data may not fully capture the solution behavior beyond it. 
Nevertheless, the EFR-OpInf errors remain moderate and bounded. 
In contrast, the weakly regularized OpInf error grows rapidly around $t \approx 3042$, while the strongly regularized OpInf error remains bounded through validation but diverges near $t \approx 3082$ in the testing region.
All three EFR-OpInf models remain bounded over the full interval $[3000,3100]$, with testing errors ranging from roughly $9.5\%$ to $14.6\%$. 
Among the EFR-OpInf models, EFR-OpInf-Proj and EFR-OpInf-HDF show similar relative $L^2$ error behavior, while EFR-OpInf-DF gives a slightly smaller relative $L^2$ error. 
This smaller error does not necessarily indicate better dynamical accuracy, as shown in Figure~\ref{fig:EFR_OpInf_KE_compare_mfu}, {where} its kinetic energy eventually settles to a nearly constant value rather than tracking the FOM dynamics.

\begin{figure}[htt]
  \centering
    \centering
    \includegraphics[width=\linewidth]{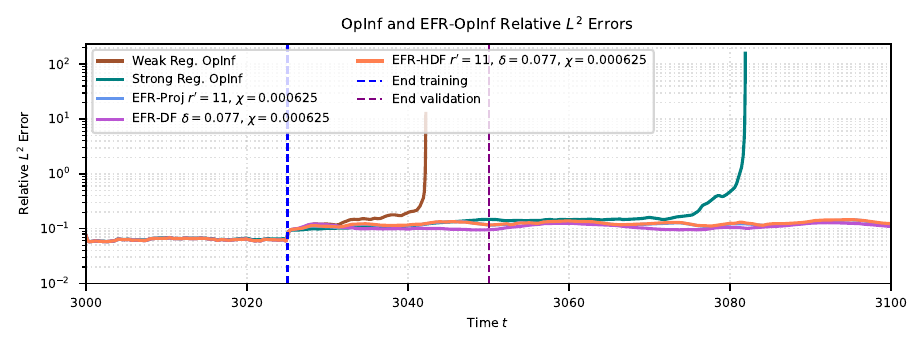}
    \caption{3D MCF at $\rm Re=5000$: Relative velocity $L^2$ errors of the OpInf and EFR-OpInf models with respect to the FOM.}
    \label{fig:EFR_OpInf_L2error_compare_mfu}
\end{figure}

Figures~\ref{fig:opinf_velocity_fields_t_3042_all} and \ref{fig:opinf_velocity_fields_t_3081_all} compare representative velocity magnitude fields in the validation and testing regions, respectively. 
At $t=3042$, the weakly regularized OpInf model has already broken down, whereas the strongly regularized OpInf and EFR-OpInf models remain bounded. 
At the later testing time $t=3081$, both standard OpInf models have become non-physical, while the EFR-OpInf models continue to produce stable velocity fields. 
These qualitative comparisons are consistent with the kinetic energy and relative $L^2$ error results. 
\begin{figure}[htt]
    \centering
    \includegraphics[width=1.0\linewidth]{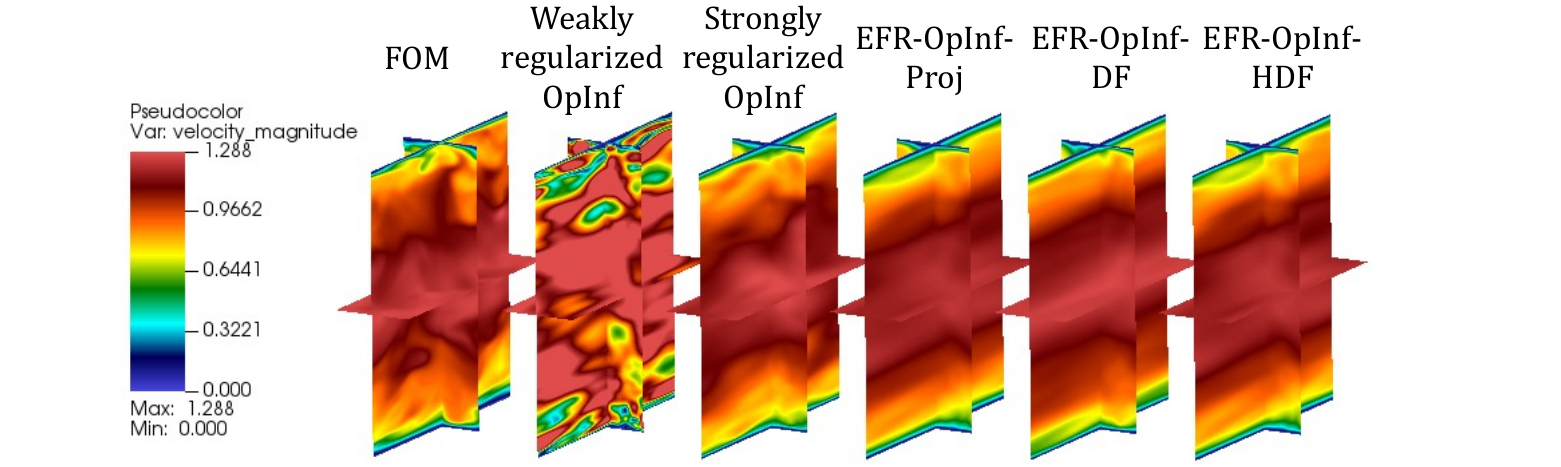}
    \caption{3D MCF at $\rm Re=5000$: Velocity magnitude fields computed by the FOM and predicted by the OpInf and EFR-OpInf models at the validation time $t=3042$.}
  \label{fig:opinf_velocity_fields_t_3042_all}
\end{figure}

\begin{figure}[htt]
    \centering
    \includegraphics[width=1.0\linewidth]{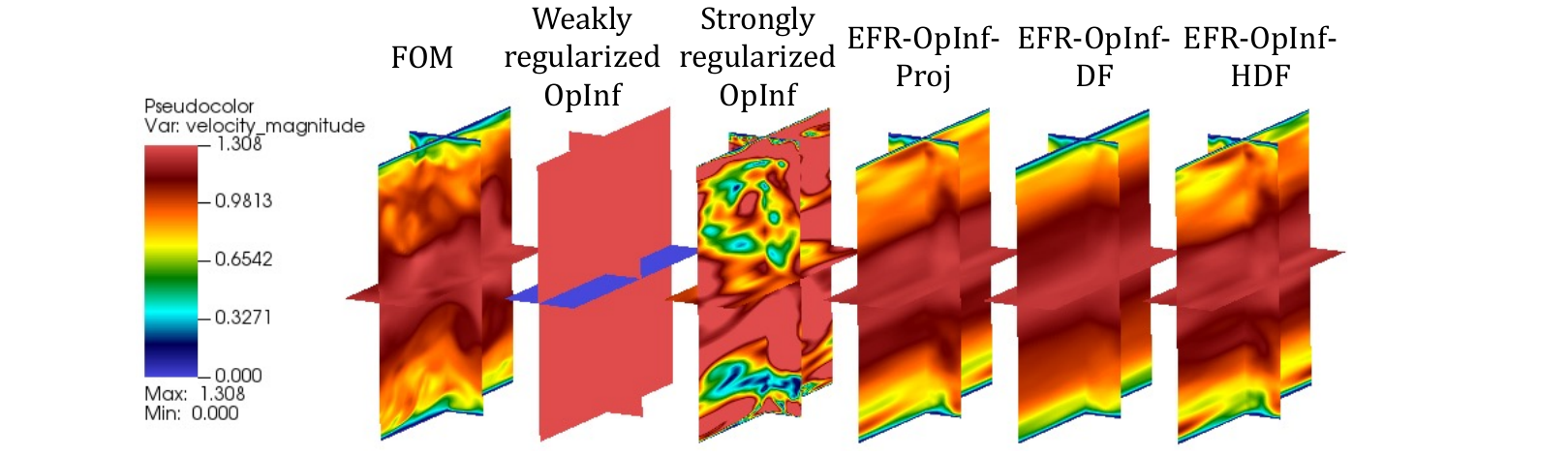}
    \caption{3D MCF at $\rm Re=5000$: Velocity magnitude fields computed by the FOM and predicted by the OpInf and EFR-OpInf models at the testing time $t=3081$.}
  \label{fig:opinf_velocity_fields_t_3081_all}
\end{figure}

Finally, Figure~\ref{fig:Pareto_Vary_Chi_MFU} compares the mean relative kinetic energy error and online ROM cost for different EFR-OpInf models obtained by varying the relaxation parameter $\chi$, while fixing the projection cutoff at $r'=11$ and the differential filter radius at $\delta=0.077$. 
Figure~\ref{fig:Pareto_Vary_Chi_MFU}~(a) reports the mean relative kinetic energy error over the training and validation regions, while Figure~\ref{fig:Pareto_Vary_Chi_MFU}~(b) includes the testing region. 
The comparison shows that some parameter choices can appear accurate over the shorter {time} interval but become less reliable once the testing region is included. 
Faded markers correspond to EFR-OpInf-DF runs that become over-diffusive, with kinetic energy settling to a nearly constant value. 
Although these runs may have small mean relative kinetic energy errors, they do not accurately reproduce the temporal variation of the FOM kinetic energy. 
Overall, EFR-OpInf-HDF performs similarly to EFR-OpInf-Proj in terms of the mean relative kinetic energy error, while avoiding the over-diffusive behavior observed for EFR-OpInf-DF.
It is also important to note that the online costs of the EFR-OpInf models are close to each other and comparable to that of the strongly regularized OpInf baseline. 
In particular, the wall-clock times for a simulation over the full time interval $[3000, 3100]$ remain within a narrow range, approximately between $6.0$ and $7.4$ seconds. 
Hence, for this scenario, the improved long-time stability provided by the EFR-OpInf models is obtained without a substantial increase in online computational cost.
Moreover, compared with the FOM, which requires approximately $981$ seconds to simulate the same time interval using nekRS on a single NVIDIA TITAN V GPU, the EFR-OpInf models achieve a wall-clock speedup of approximately $132\times$--$163\times$.

\begin{figure}[!h]
  \centering
    \begin{subfigure}[t]{0.49\linewidth}
    \centering
    \includegraphics[width=\linewidth]{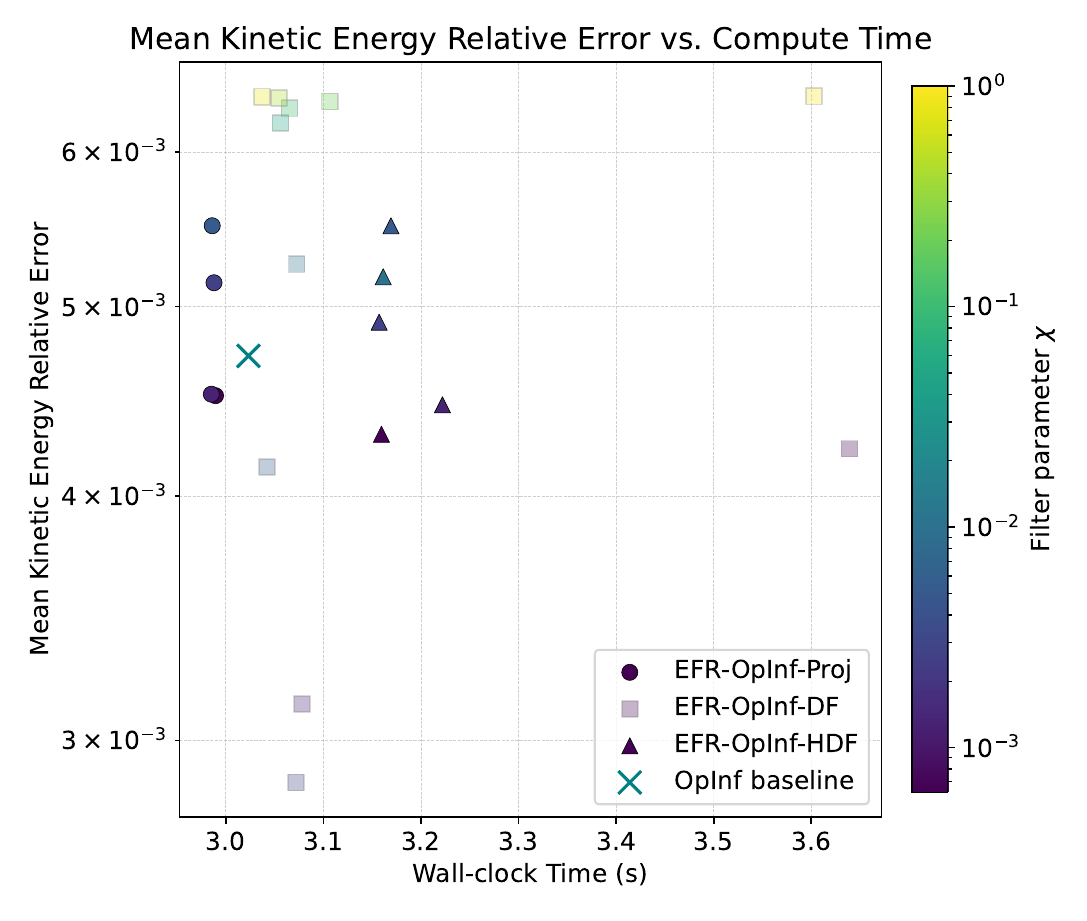}
    \caption{Training and validation regions}
  \end{subfigure}\hfill 
  \begin{subfigure}[t]{0.49\linewidth}
    \centering
    \includegraphics[width=\linewidth]{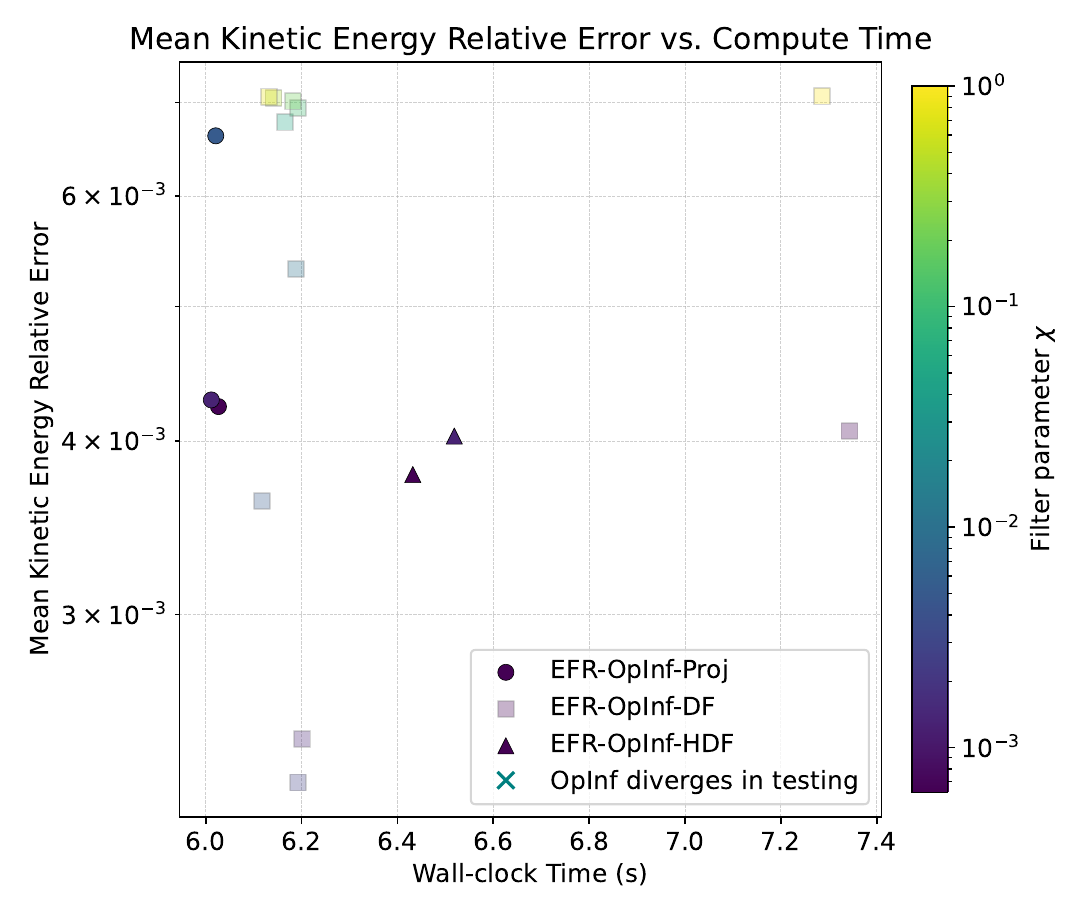}
    \caption{Training, validation, and testing regions}
  \end{subfigure}
  \caption{3D MCF at $\mathrm{Re}=5000$: Pareto comparison of mean relative kinetic energy error and online ROM cost for $r=29$. Panel~(a) uses the training and validation interval $t \in [3000,3050]$, while Panel~(b) includes the testing interval $t \in [3000,3100]$. The EFR-OpInf models are obtained by varying $\chi$, with $r'=11$ and $\delta=0.077$ fixed. Faded markers correspond to over-diffusive EFR-OpInf-DF runs.}
  \label{fig:Pareto_Vary_Chi_MFU}
\end{figure}

\newpage

\section{Conclusions} 
    \label{sec:conclusions}

We have proposed an EFR stabilization strategy for OpInf ROMs of complex fluid flow simulations.  
The main novelty of this EFR-OpInf framework is the LES-inspired online ROM spatial filtering
that significantly improves long-term stability and predictive performance of standard OpInf.
Another advantage of the new EFR-OpInf strategy is that it reduces, and in some cases even eliminates, the need for 
regularization in the OpInf problem. 
Furthermore, EFR-OpInf provides physical interpretability for its hyperparameters: 
each stabilization parameter can be associated 
with a concrete and observable modification of the spatial content of the ROM solution, enabling users to select an appropriate filter with minimal tuning and parameter choices guided by numerical analysis.
We have investigated 
EFR-OpInf for two of the most common ROM spatial filters: the completely non-intrusive ROM projection filter and the ROM differential filter, which requires additional FOM-level information. We have also studied the recently introduced hybrid ROM projection-differential filter, which provides minimal regularization by adjusting only the smallest scales of flow, preventing excessive diffusion. %
The proposed EFR-OpInf has been evaluated on a parameter-dependent, high-P\'eclet-number convection--diffusion--reaction problem, %
a transitional 2D flow past a cylinder, and a 3D turbulent minimal channel flow. 
For all these test problems, EFR-OpInf reduces the standard OpInf prediction errors by up to an order of magnitude, and remains stable over long prediction horizons in cases for which standard OpInf diverges.
Furthermore, EFR-OpInf's computational cost is comparable to that of standard OpInf.

In Table \ref{tab:conclusion_table}, we %
include a qualitative comparison of the performance of each of the considered filters on our test cases, with checkmarks (\checkmark) indicating a succesful strategy, and x-marks (\xmark) indicating an unsuccesful strategy. For the CDR, EFR-OpInf-DF was the clear winner, but both the projection and hybrid filters showed strong performance in accuracy and stability. For both of the NSE examples, the differential filter proved overly-diffusive, a phenomenon which can likely be explained by the difference in complexity of dynamics. 
Indeed, for all of the examples, the differential filter produced essentially a steady state solution, which was appropriate for the CDR problem but not for either of the NSE problems. 
Over-diffusivity has been a known problem for the differential filter~\cite{moore2026filtering}, but we stress that %
it has performed well for Galerkin ROMs of complex 3D flows previously~\cite{Wells2017EF}. For the NSE problems, we found that that both the projection and hybrid filters performed comparably, though the hybrid filter was more receptive to tuning (see Figure \ref{fig:Pareto_Vary_Chi}). 
Because of %
the similar performance of the projection and hybrid filters, we would suggest practitioners first attempt the fully non-intrusive and computationally inexpensive projection filter %
when seeking to adapt the new EFR-OpInf strategy.
}
\begin{table}[ h]
    \centering
    \begin{tabular}{c|c|c|c}
         & EFR-OpInf-Proj & EFR-OpInf-DF & EFR-OpInf-HDF \\
        CDR & \checkmark & \checkmark &  \checkmark    \\
         NSE FPC & \checkmark & \xmark & \checkmark \\
          NSE MCF & \checkmark & \xmark &  \checkmark
    \end{tabular}
    \caption{Spatial filter qualitative performance.}
    \label{tab:conclusion_table}
\end{table}

We emphasize that EFR-OpInf represents a first step in the development of online ROM spatial filtering strategies for OpInf, and opens the door to several natural extensions and follow-up directions.
First, we note that the experiments in this paper 
used only the most standard versions of EFR-OpInf's two main components: unstructured OpInf and EFR.
Probably the most natural direction for future work is to extend the new EFR-OpInf framework by 
incorporating more sophisticated OpInf methods 
and EFR strategies.
At the EFR level, the results presented in this paper use a simple heuristic for determining when to incorporate filtering: if the model is within the temporal training region, we assume no filtering is needed, and once the model exceeds the training regime we turn on the full amount of filtering to ensure stability. 
While this proved effective, the flexibility of the EFR-OpInf model enables much more granular control, allowing for the adjustment of filtering amounts during the online stage.
At the OpInf level, incorporating function-space OpInf formulations that consistently incorporate the discretization-induced spatial inner product~\cite{QFW21}, as well as extensions that go beyond POD-based reduction such as quadratic manifolds~\cite{GEELEN2023quadraticmanifold} constitute future research directions.
Moreover, extensions based on distributed OpInf~\cite{FARCAS2025109619} provide pathways for addressing large-scale applications.

Another direction for future work lies in ensuring long-term stability of OpInf when used in multi-component modeling. 
While long-term stability is already a desirable model property, 
it is even more critical in the deployment of surrogates as part of multi-component models such as digital twins. 
In this setting, a ROM deployed as a component of a digital twin or generic coupled model can provide significant speedups and flexibility over a traditional monolithic model~\cite{MooreCSRI,tezaur2026hybridcouplingoperatorinference}. 
OpInf is a particularly attractive ROM choice here due to its non-intrusive compatibility with legacy solvers. 
Stability is particularly critical in digital twins, because the failure of one small component can imply the failure of the entire simulation.

Yet another direction for future work is to use other LES-inspired methodologies for OpInf ROMs. 
This includes Leray type models %
\cite{Wells2017EF,xie2018-NA_Leray}, and approximate deconvolution methods~\cite{Moore2025_ADLROM,SANFILIPPO2023_ADL,XIE2017AD_ROM}. %
It is also worthwhile to raise the question of model consistency~\cite{strazzullo2021_EFRConsistency}, noting that our EFR-OpInf ROM creates an approximation for filtered coefficients while using unfiltered data. It would be interesting to see if results improve when filtered data are used, which would also provide a connection with prior OpInf filtering work in~\cite{farcas2023filtering}. 

The main contribution of this paper is the introduction and development of an online, ROM-level spatial filtering strategy for OpInf. 
This is part of a broader program for reducing dependence on 
hyperparameter tuning in data-driven modeling and scientific machine learning %
by leveraging mathematically supported and physically interpretable spatial filtering to promote stability. To our knowledge, the only other investigation in this direction occurs in~\cite{Rezaian2023_GNN_FR}, which found that spatial filtering of a differential type could improve model performance for graph neural networks modeling reactive flow. This paper also found that standard machine learning techniques, including adjusting weight decay and dropout probability, were ineffective in improving solution accuracy on their example problem. 
While more work is needed to extend EFR-OpInf to realistic flow applications, the findings in our paper and the results in~\cite{Rezaian2023_GNN_FR} 
demonstrate that exploiting spatial filtering represents a promising strategy in data-driven modeling and scientific machine learning, where 
a wealth of mathematical and physical principles developed in LES~\cite{berselli2006LES,Sagaut2006LES} are still waiting to be explored.

\bigskip

\section*{Declaration of competing interest}
The authors declare that they have no known competing financial interests or personal relationships that could have appeared
to influence the work reported in this paper.

\section*{Acknowledgments}
The %
second author was supported by the Yushan Fellow Program %
of the Ministry of Education (MOE), Taiwan (MOE-115-YSFMS-0010-001-P1).
The fourth author was supported in part by the National Science Foundation under Grant No.~DMS-2436357.
The first, third, fifth, and sixth authors %
wish to acknowledge Sandia National Laboratories' Laboratory Directed Research and Development (LDRD) program and the U.S. Department of Energy, Office of Science, Office of Advanced Scientific Computing Research, Mathematical Multifaceted Integrated Capability Centers (MMICCs) program, which was performed under Field Work Proposal 22025291 and the Multifaceted Mathematics for Predictive Digital Twins (M2dt) project.
Sandia National Laboratories is a multi-mission laboratory managed and operated by National Technology and Engineering Solutions of Sandia, LLC., a wholly owned subsidiary of Honeywell International, Inc., for the U.S. Department of Energy’s National Nuclear Security Administration under contract DE-NA0003525.

\appendix

\section{Effect of ROM Level Spatial Filtering}
    \label{appendix:spatial-filtering}

In this appendix, we discuss and demonstrate the effect of ROM level spatial filtering. 
At the FOM level, spatial filtering is at the core of LES.
ROM level spatial filtering is an interpretable framework for regularization of fluid flow. %
As discussed in Section~\ref{sec:opinf_regularization}, in the OpInf setting, Tikhonov regularization %
for a hyperparameter $\hyperparameter$ can be involved, and it is not a trivial task to determine the correct hyperparameter ranges or what the effect of choosing one hyperparameter vs. another is on the 
the performance of the resulting OpInf model. %
In contrast, the effect of ROM spatial filtering can be immediately visualized at any point in the simulation before any model tuning takes place with out of the box suggestions for filter parameter choices %
guided by numerical analysis. 

\begin{figure}[htt]
    \centering
    \includegraphics[width = 0.75\linewidth]{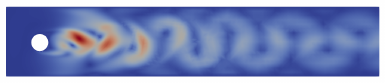}
    \caption{FPC: mean-centered FOM velocity magnitude  at $t = 8.0$.}
\label{fig:mean_centered_snapshot_t_8.0}
  \end{figure}

To demonstrate this,  for a variety of filter parameter options, we will show the effects of each of our three ROM level spatial filtering options (detailed in %
Section~\ref{subsec:filtering_methods}) on the FOM data (which is first projected onto the POD basis %
to allow easy ROM level computation) contained in the snapshot in Figure \ref{fig:Cylinder_FOM_sol}. %
We emphasize that these filtering results are not the result of any time-stepping ROM scheme; they are solely the result of applying filtering at a ROM level to a particular FOM snapshot. To make the filter effects as clear as possible, and because the actual filtering we will use is applied to mean-centered snapshots, all of the results in this section will be applied to the mean-centered snapshot of the FOM velocity field at $t = 8.0$ depicted in Figure \ref{fig:mean_centered_snapshot_t_8.0}. %

\paragraph{ROM Projection Filter}
We begin with the projection filter, which is %
closely related to the %
ROM expansion illustrated in Figure \ref{fig:rom-basis-funcs}.
A key preliminary step in projection filtering is to examine the POD basis and assess how the retained modes influence the reduced-order representation. Figure~\ref{fig:rom-basis-funcs} indicates that by approximately the \(15^\text{th}\) mode, the basis functions primarily capture increasingly fine-scale, low-energy features of the fluctuating velocity field. Retaining these higher-index modes can improve the expressiveness and in-sample accuracy of the OpInf ROM; however, these same components can promote numerical instability in out-of-sample predictions. The main advantage of EFR-OpInf-Proj is 
its ability to use these small-scale components to produce a more accurate and expressive ROM, while limiting their deleterious effects outside of the training regime in a completely non-intrusive way. 

\begin{figure}[httbp]
  \centering

  \begin{subfigure}{0.48\linewidth}
    \centering
    \includegraphics[width=\linewidth]{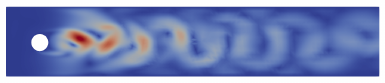}
    \caption{Effect of projection filter for $r' = 20$ (no filtering occurs).}
  \end{subfigure}

  \vspace{0.8em} %

  \begin{subfigure}{0.48\linewidth}
    \centering
    \includegraphics[width=\linewidth]{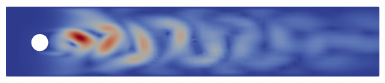}
    \caption{Effect of projection filter for $r' = 10$.}
  \end{subfigure}\hfill
  \begin{subfigure}{0.48\linewidth}
    \centering
    \includegraphics[width=\linewidth]{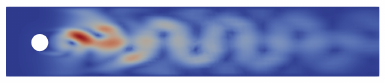}
    \caption{Effect of projection filter for $r' = 5$.}
  \end{subfigure}

  \caption{%
  FPC: effect of ROM projection filter on magnitude of fluctuating component of velocity field at time $t = 8.0$, color on identical scale.}
  \label{fig:cylinder_proj_filtering}
\end{figure}

Figure \ref{fig:cylinder_proj_filtering} demonstrates the 
effect of 
projection filtering 
on the FOM velocity magnitude at time $t = 8.0$ by selectively dropping ROM modes, and displays the same patterns as found in the POD deconstruction. When the filter cutoff is the same as the maximum POD dimension ($r' = r = 20$), the ``filtered'' solution (no filtering 
actually occurs for $r' = r)$ contains most of the same fine-structure details as the FOM velocity magnitude. The 10-dimensional projection filter diffuses the velocity representation mostly in the right half of the domain, whereas the 5-dimensional projection filter begins to merge the vortices that form directly after the cylinder. While these effects are small individually, 
they compound when incorporated into a time-stepping scheme of several hundred iterations. 
We note that it is important to be judicious with filtering: while $r' = 5$ can work for stabilization of the solution, it could also lead to excessive diffusion being added to the solution. The right approach to filtering for EFR is to determine how to add just enough to stabilize the solution. 

\begin{figure}[htt]
  \centering

  \begin{subfigure}[t]{0.48\linewidth}
    \centering
    \includegraphics[width=\linewidth]{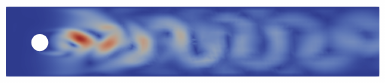}
    \caption{Effect of differential filter for $\delta = 0.01$. }
  \end{subfigure}\hfill
  \begin{subfigure}[t]{0.48\linewidth}
    \centering
    \includegraphics[width=\linewidth]{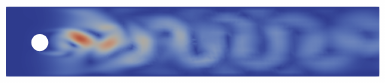}
    \caption{Effect of differential filter for $\delta = 0.02$.}
  \end{subfigure}
  \caption{%
  FPC: magnitudes of fluctuating components of filtered ROM using differential filter at $t = 8.0$ for ROM basis dimension of $20$ for $\delta = 0.01$ and $0.02$, colored on identical scale.
  }
  \label{fig:ROM_DF_Compare}
\end{figure}

\paragraph{ROM Differential Filter}
Next, we will move on to the differential filter. 
The projection and differential filters act by fundamentally different mechanisms, which explains the distinct differences between projection filtering in Figure~\ref{fig:cylinder_proj_filtering} and differential filtering in Figure~\ref{fig:ROM_DF_Compare}. The projection filter eliminates high-index modes while leaving the retained, energetically dominant modes completely intact. This leaves the velocity field's pointwise magnitude essentially unchanged in the wake of the cylinder. By contrast, the differential filter in Figure \ref{fig:ROM_DF_Compare} solves an elliptic PDE and therefore performs a spatially continuous smoothing whose local strength depends on the solution gradients; as a result it tends to reduce the magnitude of the velocity field in regions of large gradients, such as the cylinder wake. 
Figures~\ref{fig:cylinder_proj_filtering} and~\ref{fig:ROM_DF_Compare} show that the projection filter eliminates small-scale components but preserves the amplitudes of the large coherent structures, whereas the differential filter attenuates magnitudes where gradients are large throughout the entire domain. 
We note that this level of filtering throughout the domain can introduce excessive amounts of diffusion into the solution. 
Indeed, Figure \ref{fig:ROM_DF_Compare} shows post-cylinder vortices with lower pointwise velocities than for the projection filter, and this diminishing energy will have effects downstream that tend towards over-diffusion (see \cite{Moore2025_ADLROM,SANFILIPPO2023_ADL} for over-diffusivity of the differential filter in a different setting). 

\begin{figure}[httbp]
  \centering

  \begin{subfigure}[t]{0.48\linewidth}
    \centering
    \includegraphics[width=\linewidth]{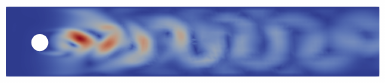}
    \caption{Effect of hybrid differential filter for $\delta = 0.01$.}
  \end{subfigure}

  \vspace{0.8em} %

  \begin{subfigure}[t]{0.48\linewidth}
    \centering
    \includegraphics[width=\linewidth]{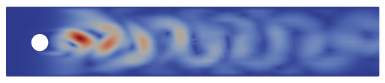}
    \caption{Effect of hybrid differential filter for $\delta = 0.02$.}
  \end{subfigure}\hfill
  \begin{subfigure}[t]{0.48\linewidth}
    \centering
    \includegraphics[width=\linewidth]{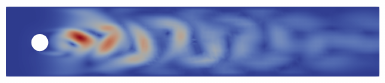}
    \caption{Effect of hybrid differential filter for $\delta = 0.2$.}
  \end{subfigure}

  \caption{%
  FPC: effect of hybrid differential filter on magnitude of FOM velocity at $t = 8.0$ for ROM basis dimension of $20$ for $\delta = 0.01$, $0.02$, and $0.2$, colored on identical scale.}
  \label{fig:ROM_DF_partial_Compare}
\end{figure}

\paragraph{Hybrid ROM Projection-Differential Filter}
To address the potential over-diffusivity of the differential filter, we use the recently proposed hybrid ROM projection-differential filter~\cite{strazzullo2025variational}.
Consider the application of the hybrid projection-differential filter to the FOM velocity magnitude at $t = 8.0$, plotted in Figure \ref{fig:ROM_DF_partial_Compare}. 
We note that there is very little difference in the performance of the filter for $\delta = 0.01$ and $\delta = 0.02$, and there is no visible diffusion in the post-wake structures. This is notably different from the results of the differential filter for the same values of $\delta$ in Figure \ref{fig:ROM_DF_Compare}. Even for the much larger value of $\delta = 0.2$, almost no diffusion 
appears in the wake for the hybrid projection-differential filter, 
suggesting that this filter could correct the over-diffusivity observed in the EFR-OpInf-Proj and EFR-OpInf-DF models in Section~\ref{sec:numerical_results}. 

\section*{Code / Data Availability}
The code used to produce the CDR and NSE results is available at \href{https://github.com/irmoore/EFR_Opinf}{https://github.com/irmoore/EFR\_Opinf}, along with instructions for ------------------------reproducing the results. This code is based on the \texttt{opinf}~\cite{opinf_python} python package and the FEniCSx~\cite{BarattaEtal2023} project.

\bibliographystyle{siam}
\bibliography{EFR}

\end{document}

%% file: Notation.tex
\newcommand{\peclet}{P\'{e}clet }

\newcommand{\umean}{\mathbf{u}_{\text{mean}}}
\newcommand{\ubar}{\bar{\mathbf{u}}}
\newcommand{\basisfunc}{\boldsymbol{\varphi}^r}
\newcommand{\basisvec}{\boldsymbol{\phi}^r}
\newcommand{\PODbasis}{\boldsymbol{\Phi}}
\newcommand{\romfunc}{\mathbf{u}_r}
\newcommand{\filfunc}{\overline{\mathbf{u}}_r}

\newcommand{\vfunc}{\mathbf{v}_{r'}}
\newcommand{\vcoef}{\hat{v}}
\newcommand{\vstate}{\hat{\mathbf{v}}}

\newcommand{\wfunc}{\mathbf{w}_{r}}

\newcommand{\state}{\hat{\mathbf{u}}}
\newcommand{\oneDstate}{\hat{u}}
\newcommand{\statecoef}{\hat{u}}
\newcommand{\filstate}{\overline{\hat{\mathbf{u}}}}
\newcommand{\filcoef}{\overline{\hat{u}}}

\newcommand{\linearOP}{\widehat{A}}
\newcommand{\quadOP}{\widehat{H}}
\newcommand{\constantOP}{\hat{\mathbf{c}}}

\newcommand{\GROMlinearOP}{A}
\newcommand{\GROMquadOP}{H}
\newcommand{\GROMconstantOP}{\mathbf{c}}

\newcommand{\hyperparameter}{\boldsymbol{\beta}}

\newcommand{\bu}{\mathbf{u}}
\newcommand{\bx}{\mathbf{x}}
\newcommand{\obu}{\overline{\mathbf{u}}}

\newcommand{\ddt
}[1]{\frac{\partial {#1}}{\partial t}}

\newtheorem{remark}{Remark}
\newtheorem{theorem}{Theorem}[section]
\newtheorem{definition}{Definition}[section]
\newtheorem{corollary}{Corollary}[theorem]
\newtheorem{lemma}[theorem]{Lemma}

%% file: EFR.bib
@article{ahmed2021closures,
  title={On closures for reduced order models $-$ A spectrum of first-principle to machine-learned avenues},
  author={Ahmed, S. E. and Pawar, S. and San, O. and Rasheed, A. and Iliescu, T. and Noack, B. R.},
  journal={Phys. Fluids},
  fjournal={Physics of Fluids},
  volume={33},
  number={9},
  pages={091301},
  year={2021},
  publisher={AIP Publishing LLC}
}

@book{rozza2022advanced,
  title={Advanced reduced order methods and applications in computational fluid dynamics},
  author={Rozza, G. and Stabile, G. and Ballarin, F.},
  year={2022},
  publisher={SIAM}
}

@book{noack2011reduced,
  title={Reduced-Order Modelling for Flow Control},
  author={Noack, B. R. and Morzynski, M. and Tadmor, G.},
  volume={528},
  year={2011},
  publisher={Springer Verlag}
}

@book{hesthaven2015certified,
  title={Certified Reduced Basis Methods for Parametrized Partial Differential Equations},
  author={Hesthaven, J. S. and Rozza, G. and Stamm, B.},
  journal={SpringerBriefs in Mathematics},
  year={2015},
  publisher={Springer}
}

@book{quarteroni2015reduced,
  title={Reduced Basis Methods for Partial Differential Equations: An Introduction},
  author={Quarteroni, A. and Manzoni, A. and Negri, F.},
  volume={92},
  year={2015},
  publisher={Springer}
}

@PhdThesis{moore2026filtering,
  author = 	 {Moore, I.},
  title = 	 {Filtering and Domain Decomposition Techniques for Intrusive and Non-intrusive Reduced Order Models of Convection-Dominated Problems},
  school = 	 {Virginia Tech},
  year = 	 {2026},
  OPTkey = 	 {},
  OPTtype = 	 {},
  OPTaddress = 	 {},
  OPTmonth = 	 {},
  OPTnote = 	 {},
  OPTannote = 	 {}
}

@Book{holmes2012turbulence,
  author = 	 {Holmes, P. and Lumley, J. L. and Berkooz, G. and Rowley, C. W.},
  title = 	 {Turbulence, Coherent Structures, Dynamical Systems and Symmetry},
  publisher = 	 {Cambridge},
  year = 	 {2012},
  OPTkey = 	 {},
  OPTvolume = 	 {},
  OPTnumber = 	 {},
  OPTseries = 	 {},
  OPTaddress = 	 {},
  edition = 	 {second},
  OPTmonth = 	 {},
  OPTnote = 	 {},
  OPTannote = 	 {}
}

@article{strazzullo2025variational,
  title={Variational multiscale evolve and filter strategies for convection-dominated flows},
  author={Strazzullo, M. and Ballarin, F. and Iliescu, T. and Rebollo, T. C.},
  journal={Comput. Meth. Appl. Mech. Eng.},
  fulljournal={Computer methods in applied mechanics and engineering},
  volume={438},
  number={11781},
  pages={117811},
  year={2025}
}

@article{volkwein2013proper,
  title={Proper orthogonal decomposition: Theory and reduced-order modelling},
  author={Volkwein, S.},
  journal={Lecture Notes, University of Konstanz},
  url={http://www.math.uni-konstanz.de/numerik/personen/volkwein/teaching/POD-Book.pdf},
  year={2013}
}

@misc{vijaywargiya2025tensorparametrichamiltonianoperator,
      title={Tensor parametric {H}amiltonian operator inference}, 
      author={Arjun Vijaywargiya and Shane A. McQuarrie and Anthony Gruber},
      year={2025},
      eprint={2502.10888},
      archivePrefix={arXiv},
      primaryClass={math.NA},
      url={https://arxiv.org/abs/2502.10888}, 
}

@article{GEELEN2023quadraticmanifold,
title = {Operator inference for non-intrusive model reduction with quadratic manifolds},
journal = {Computer Methods in Applied Mechanics and Engineering},
volume = {403},
pages = {115717},
year = {2023},
issn = {0045-7825},
doi = {https://doi.org/10.1016/j.cma.2022.115717},
url = {https://www.sciencedirect.com/science/article/pii/S0045782522006727},
author = {Rudy Geelen and Stephen Wright and Karen Willcox}
}

@inbook{Rezaian2023_GNN_FR,
author = {Elnaz Rezaian and Karthikeyan Duraisamy},
title = {Predictive Modeling of Complex Flows using Regularized Conditionally Parameterized Graph Neural Networks},
booktitle = {AIAA SCITECH 2023 Forum},
year = {2023},
chapter = {},
pages = {},
doi = {10.2514/6.2023-1284},
URL = {https://arc.aiaa.org/doi/abs/10.2514/6.2023-1284}
}

@article{Mohebujjaman2019physically,
author = {Mohebujjaman, M. and Rebholz, L.G. and Iliescu, T.},
title = {Physically constrained data-driven correction for reduced-order modeling of fluid flows},
journal = {International Journal for Numerical Methods in Fluids},
volume = {89},
number = {3},
pages = {103-122},
doi = {https://doi.org/10.1002/fld.4684},
url = {https://onlinelibrary.wiley.com/doi/abs/10.1002/fld.4684},
year = {2019}
}

@article{gunzburger2019EFR_ROM,
author = {Gunzburger, M. and Iliescu, T. and Mohebujjaman, M. and Schneier, M.},
title = {An Evolve-Filter-Relax Stabilized Reduced Order Stochastic Collocation Method for the Time-Dependent {N}avier--{S}tokes Equations},
journal = {SIAM/ASA Journal on Uncertainty Quantification},
volume = {7},
number = {4},
pages = {1162-1184},
year = {2019},
doi = {10.1137/18M1221618},

URL = { 
    
        https://doi.org/10.1137/18M1221618
    
    

}
}

@article{BALAJEWICZ2016224,
title = {Minimal subspace rotation on the {S}tiefel manifold for stabilization and enhancement of projection-based reduced order models for the compressible {N}avier–{S}tokes equations},
journal = {Journal of Computational Physics},
volume = {321},
pages = {224-241},
year = {2016},
issn = {0021-9991},
doi = {https://doi.org/10.1016/j.jcp.2016.05.037},
url = {https://www.sciencedirect.com/science/article/pii/S0021999116301826},
author = {Maciej Balajewicz and Irina Tezaur and Earl Dowell}
}

@article{reyes2025timerelaxationNA,
author = {Reyes, Jorge and Tsai, Ping-Hsuan and Novo, Julia and Iliescu, Traian},
title = {A Priori Error Bounds and Parameter Scalings for the Time-Relaxation Reduced Order Model},
journal = {SIAM Journal on Numerical Analysis},
volume = {63},
number = {5},
pages = {1933-1961},
year = {2025},
doi = {10.1137/24M1708656},

URL = { 
    
        https://doi.org/10.1137/24M1708656
    
    

}
}

@article{ISSAN2023shiftedopinf,
title = {Predicting solar wind streams from the inner-heliosphere to Earth via shifted operator inference},
journal = {Journal of Computational Physics},
volume = {473},
pages = {111689},
year = {2023},
issn = {0021-9991},
doi = {https://doi.org/10.1016/j.jcp.2022.111689},
url = {https://www.sciencedirect.com/science/article/pii/S0021999122007525},
author = {Opal Issan and Boris Kramer}
}

@inbook{Qian2019TranformandLearn,
author = {Elizabeth Qian and Boris Kramer and Alexandre N. Marques and Karen E. Willcox},
title = {Transform \& Learn: A data-driven approach to nonlinear model reduction},
booktitle = {AIAA Aviation 2019 Forum},
chapter = {},
pages = {},
doi = {10.2514/6.2019-3707},
URL = {https://arc.aiaa.org/doi/abs/10.2514/6.2019-3707},
}

@article{QIAN20201LiftandLearn,
title = {Lift \& Learn: Physics-informed machine learning for large-scale nonlinear dynamical systems},
journal = {Physica D: Nonlinear Phenomena},
volume = {406},
pages = {132401},
year = {2020},
issn = {0167-2789},
doi = {https://doi.org/10.1016/j.physd.2020.132401},
url = {https://www.sciencedirect.com/science/article/pii/S0167278919307651},
author = {Elizabeth Qian and Boris Kramer and Benjamin Peherstorfer and Karen Willcox}
}

@article{mcquarrie2023parametric,
author = {McQuarrie, Shane A. and Khodabakhshi, Parisa and Willcox, Karen E.},
title = {Nonintrusive Reduced-Order Models for Parametric Partial Differential Equations via Data-Driven Operator Inference},
journal = {SIAM Journal on Scientific Computing},
volume = {45},
number = {4},
pages = {A1917-A1946},
year = {2023},
doi = {10.1137/21M1452810},

URL = { 
    
        https://doi.org/10.1137/21M1452810
    
    

}
}

@article{Benner2015ROM,
author = {Benner, Peter and Gugercin, Serkan and Willcox, Karen},
title = {A Survey of Projection-Based Model Reduction Methods for Parametric Dynamical Systems},
journal = {SIAM Review},
volume = {57},
number = {4},
pages = {483-531},
year = {2015},
doi = {10.1137/130932715},
URL = { 
    https://doi.org/10.1137/130932715
}
}

@misc{opinf_python,
  author       = {{Operator Inference Project Contributors}},
  title        = {{OpInf: Operator Inference in Python}},
  year         = {2026},
  howpublished = {\url{https://github.com/operator-inference/opinf}},
  note         = {Python package for data-driven, non-intrusive reduced-order modeling using Operator Inference},
}

@article{Kalashnikova:2014,
title = {Construction of energy-stable projection-based reduced order models},
journal = {Applied Mathematics and Computation},
volume = {249},
pages = {569-596},
year = {2014},
issn = {0096-3003},
doi = {https://doi.org/10.1016/j.amc.2014.10.073},
url = {https://www.sciencedirect.com/science/article/pii/S0096300314014489},
author = {Irina Kalashnikova and Matthew F. Barone and Srinivasan Arunajatesan and Bart G. {van Bloemen Waanders}}
}

@book{Benner2021ROM,
editor = {Benner, Peter and Schilders, Wil and Grivet-Talocia, Stefano and Quarteroni, Alfio and Rozza, Gianluigi and Miguel Silveira, Lu{\'\i}s},
url = {https://doi.org/10.1515/9783110671490},
title = {Model Order Reduction Volume 2: Snapshot-Based Methods and Algorithms},
publisher = {De Gruyter},
address = {Berlin, Boston},
doi = {doi:10.1515/9783110671490},
isbn = {9783110671490},
year = {2021},
lastchecked = {2026-04-10}
}

@article{Gunzburger:2007,
title = {Reduced-order modeling of time-dependent {PDE}s with multiple parameters in the boundary data},
journal = {Computer Methods in Applied Mechanics and Engineering},
volume = {196},
number = {4},
pages = {1030-1047},
year = {2007},
issn = {0045-7825},
doi = {https://doi.org/10.1016/j.cma.2006.08.004},
url = {https://www.sciencedirect.com/science/article/pii/S0045782506002337},
author = {Max D. Gunzburger and Janet S. Peterson and John N. Shadid},
}

@Article{Parish2025residual,
author={Parish, Eric
and Yano, Masayuki
and Tezaur, Irina
and Iliescu, Traian},
title={Residual-Based Stabilized Reduced-Order Models of the Transient Convection--Diffusion--Reaction Equation Obtained Through Discrete and Continuous Projection},
journal={Archives of Computational Methods in Engineering},
year={2025},
month={Apr},
day={01},
volume={32},
number={3},
pages={1885-1929},
issn={1886-1784},
doi={10.1007/s11831-024-10197-1},
url={https://doi.org/10.1007/s11831-024-10197-1}
}

@misc{Aretz:2025,
      title={Nested Operator Inference for Adaptive Data-Driven Learning of Reduced-order Models}, 
      author={Nicole Aretz and Karen Willcox},
      year={2025},
      eprint={2508.11542},
      archivePrefix={arXiv},
      primaryClass={cs.LG},
      url={https://arxiv.org/abs/2508.11542}, 
}

@article{Wells2017EF,
author = {Wells, D. and Wang, Z. and Xie, X. and Iliescu, T.},
title = {An evolve-then-filter regularized reduced order model for convection-dominated flows},
journal = {International Journal for Numerical Methods in Fluids},
volume = {84},
number = {10},
pages = {598-615},
doi = {https://doi.org/10.1002/fld.4363},
url = {https://onlinelibrary.wiley.com/doi/abs/10.1002/fld.4363},
year = {2017}
}

@article{WANG2012PODClosure,
title = {Proper orthogonal decomposition closure models for turbulent flows: A numerical comparison},
journal = {Computer Methods in Applied Mechanics and Engineering},
volume = {237-240},
pages = {10-26},
year = {2012},
issn = {0045-7825},
doi = {https://doi.org/10.1016/j.cma.2012.04.015},
url = {https://www.sciencedirect.com/science/article/pii/S0045782512001429},
author = {Zhu Wang and Imran Akhtar and Jeff Borggaard and Traian Iliescu}
}

@Article{Koc2022Verifiability,
author={Koc, Birgul
and Mou, Changhong
and Liu, Honghu
and Wang, Zhu
and Rozza, Gianluigi
and Iliescu, Traian},
title={Verifiability of the Data-Driven Variational Multiscale Reduced Order Model},
journal={Journal of Scientific Computing},
year={2022},
month={Oct},
day={08},
volume={93},
number={2},
pages={54},
issn={1573-7691},
doi={10.1007/s10915-022-02019-y},
url={https://doi.org/10.1007/s10915-022-02019-y}
}

@Article{Reyes2025Verifiability,
author={Reyes, Jorge
and Tsai, Ping-Hsuan
and Moore, Ian
and Liu, Honghu
and Iliescu, Traian},
title={Verifiability and Limit Consistency of Eddy Viscosity Large Eddy Simulation Reduced Order Models},
journal={Journal of Scientific Computing},
year={2025},
month={Nov},
day={04},
volume={105},
number={3},
pages={78},
issn={1573-7691},
doi={10.1007/s10915-025-03106-6},
url={https://doi.org/10.1007/s10915-025-03106-6}
}

@article{McQuarrie2021Regularized,
author = {Shane A. McQuarrie and Cheng Huang and Karen E. Willcox},
title = {Data-driven reduced-order models via regularised Operator Inference for a single-injector combustion process},
journal = {Journal of the Royal Society of New Zealand},
volume = {51},
number = {2},
pages = {194--211},
year = {2021},
publisher = {Taylor \& Francis},
doi = {10.1080/03036758.2020.1863237},
URL = { 
https://doi.org/10.1080/03036758.2020.1863237
}
}

@misc{tezaur2026hybridcouplingoperatorinference,
      title={Hybrid coupling with operator inference and the overlapping {S}chwarz alternating method}, 
      author={Irina Tezaur and Eric Parish and Anthony Gruber and Ian Moore and Christopher Wentland and Alejandro Mota},
      year={2026},
      eprint={2511.20687},
      archivePrefix={arXiv},
      primaryClass={math.NA},
      url={https://arxiv.org/abs/2511.20687}, 
}

@article{Mullen1999Filtering,
author = {Mullen, Julie S. and Fischer, Paul F.},
title = {Filtering techniques for complex geometry fluid flows},
journal = {Communications in Numerical Methods in Engineering},
volume = {15},
number = {1},
pages = {9-18},
doi = {https://doi.org/10.1002/(SICI)1099-0887(199901)15:1<9::AID-CNM219>3.0.CO;2-Y},
url = {https://onlinelibrary.wiley.com/doi/abs/10.1002/%28SICI%291099-0887%28199901%2915%3A1%3C9%3A%3AAID-CNM219%3E3.0.CO%3B2-Y},
year = {1999}
}

@article{FARCAS2025109619,
title = {Distributed computing for physics-based data-driven reduced modeling at scale: Application to a rotating detonation rocket engine},
journal = {Computer Physics Communications},
volume = {313},
pages = {109619},
year = {2025},
issn = {0010-4655},
doi = {https://doi.org/10.1016/j.cpc.2025.109619},
url = {https://www.sciencedirect.com/science/article/pii/S0010465525001213},
author = {Ionuţ-Gabriel Farcaş and Rayomand P. Gundevia and Ramakanth Munipalli and Karen E. Willcox}
}

@inbook{farcas2023parametric,
author = {Ionut Farcas and Rayomand Gundevia and Ramakanth Munipalli and Karen E. Willcox},
title = {Parametric non-intrusive reduced-order models via operator inference for large-scale rotating detonation engine simulations},
booktitle = {AIAA SCITECH 2023 Forum},
chapter = {},
pages = {},
doi = {10.2514/6.2023-0172},
URL = {https://arc.aiaa.org/doi/abs/10.2514/6.2023-0172}
}

@article{XIE2017AD_ROM,
title = {Approximate deconvolution reduced order modeling},
journal = {Computer Methods in Applied Mechanics and Engineering},
volume = {313},
pages = {512-534},
year = {2017},
issn = {0045-7825},
doi = {https://doi.org/10.1016/j.cma.2016.10.005},
url = {https://www.sciencedirect.com/science/article/pii/S0045782516302481},
author = {Xuping Xie and David Wells and Zhu Wang and Traian Iliescu}
}

@inbook{farcas2023filtering,
author = {Ionut Farcas and Ramakanth Munipalli and Karen E. Willcox},
title = {On filtering in non-intrusive data-driven reduced-order modeling},
booktitle = {AIAA AVIATION 2022 Forum},
chapter = {},
pages = {},
doi = {10.2514/6.2022-3487},
URL = {https://arc.aiaa.org/doi/abs/10.2514/6.2022-3487}
}

@article{GKIMISKIS2025localizedopinf,
title = {Non-intrusive reduced-order modeling for dynamical systems with spatially localized features},
journal = {Computer Methods in Applied Mechanics and Engineering},
volume = {444},
pages = {118115},
year = {2025},
issn = {0045-7825},
doi = {https://doi.org/10.1016/j.cma.2025.118115},
url = {https://www.sciencedirect.com/science/article/pii/S0045782525003871},
author = {Leonidas Gkimisis and Nicole Aretz and Marco Tezzele and Thomas Richter and Peter Benner and Karen E. Willcox}
}

@misc{mcquarrie2025activelearningdatadrivenreduced,
      title={Active learning for data-driven reduced models of parametric differential systems with {B}ayesian operator inference}, 
      author={Shane A. McQuarrie and Mengwu Guo and Anirban Chaudhuri},
      year={2025},
      eprint={2601.00038},
      archivePrefix={arXiv},
      primaryClass={stat.ML},
      url={https://arxiv.org/abs/2601.00038}, 
}

@inproceedings{MooreCSRI,
    author = {Moore, I. and Wentland, C. R. and Gruber, A. and Tezaur, I.},
    title = {Domain Decomposition-Based Coupling of Operator Inference Reduced Order Models via the {S}chwarz Alternating Method},
    booktitle = {Computer Science Research Institute Summer Proceedings},
    year = {2024},
    pages = {109-127},
    publisher = {Sandia National Laboratories},
    note = {{T}echnical {r}eport SAND2024-16688O},
    editor = {Adams, M.B.P and Casey, T.A. and Reuter, B.W.}
}

@misc{DokkenFPC,
title = {Test problem 2: Flow past a cylinder ({DFG} 2{D}-3 benchmark)},
author = {Dokken, Jørgen S.},
url = {https://jsdokken.com/dolfinx-tutorial/chapter2/ns_code2.html},
note = {Accessed 04/13/2026},
year = {2026}
}

@misc{BarattaEtal2023,
  title     = {{DOLFINx}: the next generation {FEniCS} problem solving environment},
  author    = {Baratta, Igor A. and Dean, Joseph P. and Dokken, J{\o}rgen S. and Habera, Michal and Hale, Jack S. and Richardson, Chris N. and Rognes, Marie E. and Scroggs, Matthew W. and Sime, Nathan and Wells, Garth N.},
  doi       = {10.5281/zenodo.10447666},
  year      = {2023},
  howpublished = {preprint}
}

@article{ScroggsEtal2022,
  title     = {Construction of arbitrary order finite element degree-of-freedom maps on polygonal and polyhedral cell meshes},
  author    = {Scroggs, Matthew W. and Dokken, J{\o}rgen S. and Richardson, Chris N. and Wells, Garth N.},
  journal   = {ACM Transactions on Mathematical Software},
  year      = {2022},
  volume    = {48},
  number    = {2},
  doi       = {10.1145/3524456},
  pages     = {{18:1--18:23}},
}

@article{BasixJoss,
  title     = {Basix: a runtime finite element basis evaluation library},
  author    = {Scroggs, Matthew W. and Baratta, Igor A. and Richardson, Chris N. and Wells, Garth N.},
  journal   = {Journal of Open Source Software},
  year      = {2022},
  volume    = {7},
  number    = {73},
  doi       = {10.21105/joss.03982},
  pages     = {3982}
}

@article{GUO2022bayesianopinf,
title = {Bayesian operator inference for data-driven reduced-order modeling},
journal = {Computer Methods in Applied Mechanics and Engineering},
volume = {402},
pages = {115336},
year = {2022},
note = {A Special Issue in Honor of the Lifetime Achievements of J. Tinsley Oden},
issn = {0045-7825},
doi = {https://doi.org/10.1016/j.cma.2022.115336},
url = {https://www.sciencedirect.com/science/article/pii/S0045782522004273},
author = {Mengwu Guo and Shane A. McQuarrie and Karen E. Willcox}
}

@article{SHARMA2022Hamiltonian,
title = {Hamiltonian operator inference: Physics-preserving learning of reduced-order models for canonical {H}amiltonian systems},
journal = {Physica D: Nonlinear Phenomena},
volume = {431},
pages = {133122},
year = {2022},
issn = {0167-2789},
doi = {https://doi.org/10.1016/j.physd.2021.133122},
url = {https://www.sciencedirect.com/science/article/pii/S0167278921002682},
author = {Harsh Sharma and Zhu Wang and Boris Kramer}
}

@article{Gruber2025Hamiltonian,
author = {Gruber, Anthony and Tezaur, Irina},
title = {Variationally Consistent {H}amiltonian Model Reduction},
journal = {SIAM Journal on Applied Dynamical Systems},
volume = {24},
number = {1},
pages = {376-414},
year = {2025},
doi = {10.1137/24M1652490},
URL = { https://doi.org/10.1137/24M1652490
}
}

@article{strazzullo2021_EFRConsistency,
author = {Strazzullo, Maria and Girfoglio, Michele and Ballarin, Francesco and Iliescu, Traian and Rozza, Gianluigi},
title = {Consistency of the full and reduced order models for evolve-filter-relax regularization of convection-dominated, marginally-resolved flows},
journal = {International Journal for Numerical Methods in Engineering},
volume = {123},
number = {14},
pages = {3148-3178},
doi = {https://doi.org/10.1002/nme.6942},
url = {https://onlinelibrary.wiley.com/doi/abs/10.1002/nme.6942},
year = {2022}
}

@article{IVAGNES2026_EFR,
title = {A new data-driven energy-stable evolve-filter-relax model for turbulent flow simulation},
journal = {Computer Methods in Applied Mechanics and Engineering},
volume = {450},
pages = {118654},
year = {2026},
issn = {0045-7825},
doi = {https://doi.org/10.1016/j.cma.2025.118654},
url = {https://www.sciencedirect.com/science/article/pii/S0045782525009260},
author = {Anna Ivagnes and Toby {van Gastelen} and Syver Døving Agdestein and Benjamin Sanderse and Giovanni Stabile and Gianluigi Rozza}
}

@article{Bertagna2016_EFR,
author = {Bertagna, L. and Quaini, A. and Veneziani, A.},
title = {Deconvolution-based nonlinear filtering for incompressible flows at moderately large {R}eynolds numbers},
journal = {International Journal for Numerical Methods in Fluids},
volume = {81},
number = {8},
pages = {463-488},
doi = {https://doi.org/10.1002/fld.4192},
url = {https://onlinelibrary.wiley.com/doi/abs/10.1002/fld.4192},
year = {2016}
}

@article{KANEKO2020POD_ROM,
title = {Towards model order reduction for fluid-thermal analysis},
journal = {Nuclear Engineering and Design},
volume = {370},
pages = {110866},
year = {2020},
issn = {0029-5493},
doi = {https://doi.org/10.1016/j.nucengdes.2020.110866},
url = {https://www.sciencedirect.com/science/article/pii/S0029549320303605},
author = {Kento Kaneko and Ping-Hsuan Tsai and Paul Fischer}
}

@ARTICLE{kaneko2022LerayROM,
    
AUTHOR={Kaneko, Kento  and Fischer, Paul },
           
TITLE={Augmented reduced order models for turbulence},
          
JOURNAL={Frontiers in Physics},
          
VOLUME={Volume 10 - 2022},
  
YEAR={2022},
  
URL={https://www.frontiersin.org/journals/physics/articles/10.3389/fphy.2022.905392},
  
DOI={10.3389/fphy.2022.905392},
  
ISSN={2296-424X}}

@article{TSAI2025_timerelaxation,
title = {A time-relaxation reduced order model for the turbulent channel flow},
journal = {Journal of Computational Physics},
volume = {521},
pages = {113563},
year = {2025},
issn = {0021-9991},
doi = {https://doi.org/10.1016/j.jcp.2024.113563},
url = {https://www.sciencedirect.com/science/article/pii/S0021999124008118},
author = {Ping-Hsuan Tsai and Paul Fischer and Traian Iliescu}
}

@article{neda2012_EFR_analysis,
author = {Ervin, Vincent J. and Layton, William J. and Neda, Monika},
title = {Numerical Analysis of Filter-Based Stabilization for Evolution Equations},
journal = {SIAM Journal on Numerical Analysis},
volume = {50},
number = {5},
pages = {2307-2335},
year = {2012},
doi = {10.1137/100782048},

URL = { 
    
        https://doi.org/10.1137/100782048
}
}

@article{MOU2023_lengthscale,
title = {An energy-based lengthscale for reduced order models of turbulent flows},
journal = {Nuclear Engineering and Design},
volume = {412},
pages = {112454},
year = {2023},
issn = {0029-5493},
doi = {https://doi.org/10.1016/j.nucengdes.2023.112454},
url = {https://www.sciencedirect.com/science/article/pii/S0029549323003035},
author = {Changhong Mou and Elia Merzari and Omer San and Traian Iliescu}
}

@misc{tsai2026stabopdatadrivenstabilizationoperator,
      title={Stab{O}p: A Data-Driven Stabilization Operator for Reduced Order Modeling}, 
      author={Ping-Hsuan Tsai and Anna Ivagnes and Annalisa Quaini and Traian Iliescu and Gianluigi Rozza},
      year={2026},
      eprint={2602.07745},
      archivePrefix={arXiv},
      primaryClass={physics.flu-dyn},
      url={https://arxiv.org/abs/2602.07745}, 
}

@Inbook{Layton2012AD,
author="Layton, William J.
and Rebholz, Leo G.",
title="Approximate Deconvolution Operators and Models",
bookTitle="Approximate Deconvolution Models of Turbulence: Analysis, Phenomenology and Numerical Analysis",
year="2012",
publisher="Springer Berlin Heidelberg",
address="Berlin, Heidelberg",
pages="61--88",
isbn="978-3-642-24409-4",
doi="10.1007/978-3-642-24409-4_3",
url="https://doi.org/10.1007/978-3-642-24409-4_3"
}

@article{GIRFOGLIO2023FILTER,
title = {A hybrid projection/data-driven reduced order model for the {N}avier-{S}tokes equations with nonlinear filtering stabilization},
journal = {Journal of Computational Physics},
volume = {486},
pages = {112127},
year = {2023},
issn = {0021-9991},
doi = {https://doi.org/10.1016/j.jcp.2023.112127},
url = {https://www.sciencedirect.com/science/article/pii/S002199912300222X},
author = {Michele Girfoglio and Annalisa Quaini and Gianluigi Rozza}
}

@article{Germano1986Differential,
    author = {Germano, M.},
    title = {Differential filters of elliptic type},
    journal = {The Physics of Fluids},
    volume = {29},
    number = {6},
    pages = {1757-1758},
    year = {1986},
    month = {06},
    issn = {0031-9171},
    doi = {10.1063/1.865650},
    url = {https://doi.org/10.1063/1.865650},
}

@book{Sagaut2006LES,
author = "Sagaut, Pierre",
Title="Large Eddy Simulation for Incompressible Flows: An Introduction",
year="2006",
publisher="Springer Berlin Heidelberg",
address="Berlin, Heidelberg",
pages="15--44",
isbn="978-3-540-26403-3",
doi="10.1007/3-540-26403-5_2",
url="https://doi.org/10.1007/3-540-26403-5_2"
}

@book{berselli2006LES,
  title={Mathematics of large eddy simulation of turbulent flows},
  author={Berselli, Luigi C and Iliescu, Traian and Layton, William J},
  year={2006},
  publisher={Springer}
}

@article{xie2018-NA_Leray,
title = {Numerical analysis of the {L}eray reduced order model},
journal = {Journal of Computational and Applied Mathematics},
volume = {328},
pages = {12-29},
year = {2018},
issn = {0377-0427},
doi = {https://doi.org/10.1016/j.cam.2017.06.026},
url = {https://www.sciencedirect.com/science/article/pii/S0377042717303266},
author = {Xuping Xie and David Wells and Zhu Wang and Traian Iliescu}
}

@article{PEHERSTORFER2016_Opinf,
title = {Data-driven operator inference for nonintrusive projection-based model reduction},
journal = {Computer Methods in Applied Mechanics and Engineering},
volume = {306},
pages = {196-215},
year = {2016},
issn = {0045-7825},
doi = {https://doi.org/10.1016/j.cma.2016.03.025},
url = {https://www.sciencedirect.com/science/article/pii/S0045782516301104},
author = {Benjamin Peherstorfer and Karen Willcox}
}

@article{Kramer2024_Opinf,
   author = "Kramer, Boris and Peherstorfer, Benjamin and Willcox, Karen E.",
   title = "Learning Nonlinear Reduced Models from Data with Operator Inference", 
   journal= "Annual Review of Fluid Mechanics",
   year = "2024",
   volume = "56",
   number = "Volume 56, 2024",
   pages = "521-548",
   doi = "https://doi.org/10.1146/annurev-fluid-121021-025220",
   url = "https://www.annualreviews.org/content/journals/10.1146/annurev-fluid-121021-025220",
   publisher = "Annual Reviews",
   issn = "1545-4479",
   type = "Journal Article",
  }

@article{SANFILIPPO2023_ADL,
title = {Approximate deconvolution {L}eray reduced order model for convection-dominated flows},
journal = {Finite Elements in Analysis and Design},
volume = {226},
pages = {104021},
year = {2023},
issn = {0168-874X},
doi = {https://doi.org/10.1016/j.finel.2023.104021},
url = {https://www.sciencedirect.com/science/article/pii/S0168874X23001142},
author = {Anna Sanfilippo and Ian Moore and Francesco Ballarin and Traian Iliescu}
}

@article{Moore2025_ADLROM,
author = {Moore, Ian and Sanfilippo, Anna and Ballarin, Francesco and Iliescu, Traian},
title = {A Priori Error Bounds for the Approximate Deconvolution {L}eray Reduced Order Model},
journal = {Numerical Methods for Partial Differential Equations},
volume = {41},
number = {6},
pages = {e70044},
doi = {https://doi.org/10.1002/num.70044},
url = {https://onlinelibrary.wiley.com/doi/abs/10.1002/num.70044},
year = {2025}
}

@article{jimenez1991minimal,
  title={The minimal flow unit in near-wall turbulence},
  author={Jim{\'e}nez, J. and Moin, P.},
  journal={J. Fluid Mech.},
  fjournal={Journal of Fluid Mechanics},
  volume={225},
  pages={213--240},
  year={1991},
  publisher={Cambridge University Press}
}

@article{fischer2022nekrs,
  title={Nek{RS}, a {GPU}-accelerated spectral element {N}avier--{S}tokes solver},
  author={Fischer, P. and Kerkemeier, S. and Min, M. and Lan, Y.-H. and Phillips, M. and Rathnayake, T. and Merzari, E. and Tomboulides, A. and Karakus, A. and Chalmers, N. and others},
  journal={Parallel Comput.},
  fjournal={Parallel Computing},
  volume={114},
  pages={102982},
  year={2022},
  publisher={Elsevier}
}

@article{Sw19,
title = {{Projection-based model reduction: Formulations for physics-based machine learning}},
journal = {Computers \& Fluids},
volume = {179},
pages = {704--717},
year = {2019},
issn = {0045-7930},
doi = {https://doi.org/10.1016/j.compfluid.2018.07.021},
url = {https://www.sciencedirect.com/science/article/pii/S0045793018304250},
author = {Renee Swischuk and Laura Mainini and Benjamin Peherstorfer and Karen Willcox}
}

@article{QFW21,
author = {Qian, Elizabeth and Farca\c{s}, Ionu\c{t}-Gabriel and Willcox, Karen},
title = {{Reduced Operator Inference for Nonlinear Partial Differential Equations}},
journal = {SIAM Journal on Scientific Computing},
volume = {44},
number = {4},
pages = {A1934--A1959},
year = {2022},
doi = {10.1137/21M1393972}}

@incollection{fischer2017recent,
   title={Recent developments in spectral element simulations of moving-domain problems},
   author={Fischer, P. and Schmitt, M. and Tomboulides, A.},
   booktitle={Recent Progress and Modern Challenges in Applied Mathematics, Modeling and Computational Science},
   pages={213--244},
   year={2017},
   publisher={Springer}
 }

@misc{kaneko_nekrom,
  author       = {Kaneko, K. and
                  Tsai, P.-H. and
                  Christensen, N. and
                  Fischer, P.},
  title        = {{NekROM}},
  year         = 2026,
  publisher    = {Zenodo},
  version      = {v0.2.0},
  doi          = {10.5281/zenodo.18487297},
  url          = {https://doi.org/10.5281/zenodo.18487297},
  howpublished={\url{https://github.com/Nek5000/NekROM}},
}
